\documentclass{article}

\usepackage{amsfonts} 
\usepackage{amssymb,amsmath} 
\usepackage{graphicx,color}
\usepackage[breaklinks,colorlinks=true,linkcolor=blue,citecolor=red,backref=page]{hyperref}

\DeclareMathAlphabet{\itbf}{OML}{cmm}{b}{it}

\advance \textwidth by 0.9in
\advance \oddsidemargin by -0.65in
\advance \textheight by 1.2in
\advance \topmargin by -.5in

\def\RR{\mathbb{R}}

\def\by{{{\itbf y}}}
\def\bx{{{\itbf x}}}
\def\bX{{{\itbf X}}}

\newcommand{\bk}{{\itbf k}}

\def\qed{\hfill {\small $\Box$} \\}

\newtheorem{proposition}{Proposition} 
\newtheorem{lemma}{Lemma} 
\newtheorem{remark}{Remark} 

\begin{document}

\title{How a moving passive observer can perceive its environment~? \\The Unruh effect revisited}

\author{Mathias Fink\footnote{Institut Langevin, ESPCI and CNRS, PSL Research University, 1 rue Jussieu, 75005 Paris, France} \,
and Josselin Garnier\footnote{Centre de Math\'ematiques Appliqu\'ees, Ecole Polytechnique, 91128 Palaiseau Cedex, France. {\tt josselin.garnier@polytechnique.edu} (corresponding author).}}

\maketitle

\begin{abstract}
We consider a point-like observer that moves in a medium illuminated by noise sources with Lorentz-invariant spectrum.
We show that the autocorrelation function of the signal recorded by the observer allows it to perceive its environment.
More precisely,  we consider an observer with constant acceleration (along a Rindler trajectory)
and we exploit the recent work on the emergence of the Green's function 
from the cross correlation of signals transmitted by noise sources.
First we recover the result that the signal recorded by the observer has a constant Wigner transform, i.e. a constant local spectrum,
when the medium is homogeneous (this is the classical analogue of the Unruh effect).
We complete that result by showing that the Rindler trajectory is the only straight-line trajectory that satisfies this property.
We also show that, in the presence of an obstacle
{in the form of  an infinite perfect mirror}, 
the Wigner transform is perturbed when the observer comes into the neighborhood
of the obstacle. The perturbation makes it possible for the observer to determine its position relative to the obstacle
{once the entire trajectory has been traversed.}
\end{abstract}

\section{Introduction}

It has been shown that the cross correlations of the signals recorded by a stationary receiver array and transmitted by
opportunistic or ambient noise sources can be used to image the environment of the array \cite{garpapa09,garpapa16}.
It is possible to estimate travel times between receivers in order to estimate the background velocity tomographically \cite{shapiro05}
or to detect and localize reflectors in the medium \cite{gouedard08}. These ideas have been exploited in particular 
{in} seismology \cite{curtis06,wap10a} and they require the use of arrays or networks of sensors.
Indeed, the autocorrelation function of the signal recorded by a unique stationary point-like receiver contains little information
about its environment. We will see that a stationary receiver or observer can, however, estimate its distance from an obstacle but nothing more. 

The situation is different when the observer is moving, because the observer may then be able to exploit the synthetic aperture
generated by its trajectory.
For instance, active synthetic aperture radar (when the moving antenna transmits and receives) has proved to be a very efficient imaging modality \cite{cheney01,curlander}
and bistatic or passive versions (when the moving antenna only records and uses signals transmitted by controlled or opportunistic sources)
are now the subject of intense research \cite{antoniou12,finkgar15,garpapa15,rodriguez10}.
When the observer is moving, the autocorrelation function of the recorded signal depends in a complicated and interesting way of the 
environment. In \cite{fink17} the situation in which a receiver is moving along a circular trajectory is addressed.  
It is shown that the autocorrelation function of the recorded signal is related to
the matrix of Green's function between pairs of points along the trajectory, more
exactly  to  a  diagonal  band  of  this  matrix  whose  thickness  is  determined  by  the
velocity of the receiver.  As an application, when a point-like
reflector is present within the circular trajectory of the receiver, it is shown 
how to use the autocorrelation function of the recorded signal to localize
the reflector by migration.  The processing is, however, quite complex, because it requires
to extract the components of the correlation function due to the reflector.
It would be of great interest to determine and study a trajectory for which the unperturbed
autocorrelation function has a very simple form that makes it possible to detect and identify easily any perturbation.

By investigating the trajectories that would satisfy the desired properties, 
we have discovered connections with quantum physics and the celebrated Unruh effect.
The Unruh effect \cite{unruh} predicts that an observer along a Rindler trajectory with constant acceleration perceives the quantum vacuum as thermal radiation (i.e. it observes a thermal spectrum of particle excitations).
This effect and its applications to black hole radiation and quantum field theory
are extensively studied in the literature \cite{book:ulf,crispino}.
Moreover, it is shown to arise from the classical correlation of noise in \cite{boyer80}, using a representation of the ambient field
as a superposition of incoherent plane waves.
{In \cite{boyer80} the case of an open medium is addressed,  
but his approach could be extended to the case of a homogeneous half-space following the decomposition method introduced in \cite{rytov}.}
The same approach was used recently in \cite{ulf} 
in which the classical analog of the Unruh effect is extensively discussed and 
a simple experiment on water waves that corroborates the idea is presented.
The authors insist on the classical root of the Unruh effect as the correlation of noise in space and time. 
This work and our own research are, therefore, related, although the motivations are different.
We may say that our paper revisits the pioneering work \cite{boyer80} 
by using another approach that considers that the ambient field recorded by the observer is generated by noise sources.
This gives the same situation as in  \cite{boyer80} in a homogeneous medium,
but our approach makes it possible to consider the case where the medium is not homogeneous.
We demonstrate that the Wigner transform of the field recorded by an observer along a Rindler trajectory is constant when the medium
is homogeneous, as originally shown in \cite{boyer80} and observed in \cite{ulf}.
We  show in section \ref{sec:1} {the original result} that the Rindler trajectory is in fact the unique straight-line trajectory that satisfies this property.
We finally show that the Wigner transform of the field recorded by a Rindler observer
 is perturbed by an obstacle when the observer comes into its neighborhood.
The perturbation is analyzed in detail in section \ref{sec:2} {when the obstacle has the form of an infinite perfect mirror} and it is shown that it can be processed to extract the position of the obstacle
relative to the observer {once the entire trajectory has been traversed.}

\section{A Rindler observer in a three-dimensional open medium}
\label{sec:1}
In this section we show that an observer with a Rindler trajectory is an ideal candidate to probe the environment 
because the local spectrum of the recorded signal, in the absence of any obstacle, is independent of the position or time along the trajectory.
We also show that the Rindler trajectory is, in fact, the unique straight-line  trajectory that satisfies this property.

\subsection{The noise sources and wave fields}
We consider the three-dimensional scalar wave equation for the scalar wave field $u(t,\bx)$ in the full space $\RR^3$:
\begin{equation}
\frac{1}{c_o^2} \partial_t^2 u -\Delta u =n(t,\bx) ,
\end{equation}
with radiation condition at infinity.
The source term $n(t,\bx)$ models a noise source distribution. It is a zero-mean process, stationary in time
and delta-correlated in space:
\begin{equation}
\left< n(t,\bx)n(t',\bx') \right> = F(t-t')\delta(\bx-\bx') K(\bx) ,
\end{equation}
where $K(\bx)$ is the function that characterizes the spatial support of the noise source distribution
and  the Fourier transform $\hat{F}(\omega)$ of $F(t)$ is the power spectral density of the sources.

The analysis of the autocorrelation of the recorded signal 
follows the lines of the recent work on the emergence of the Green's function 
from the cross correlation of signals transmitted by  noise sources \cite{garpapa16}.
The covariance function of the wave field has the form
\begin{align*}
\left< u(t,\bx) u(t',\bx') \right>
=
\frac{1}{2\pi} \int_\RR \int_{\RR^3} \overline{\hat{G}\big(\omega,\bx,\by\big)}
\hat{G}\big(\omega,\bx',\by\big)
K(\by) \hat{F}(\omega) \exp\big[ i \omega ( t-t') \big]
d\by d\omega  ,
\end{align*}
where $\hat{G}(\omega,\bx,\by)$
is the three-dimensional homogeneous Green's function:
\begin{equation}
\label{eq:green3d}
\hat{G}(\omega,\bx,\by) = \frac{1}{4\pi |\bx-\by|} \exp\Big( i \frac{\omega}{c_o} |\bx-\by| \Big)  .
\end{equation}

If we assume that sources are far away, 
for instance, if the sources are at the surface of a ball with large radius,
$K(\bx) = \delta_{\partial B({\bf 0},L)} (\bx)$, then we can invoke Helmholtz-Kirchhoff identity~(\ref{hk})
and we get
\begin{align}
\nonumber
\left< u(t,\bx) u(t',\bx') \right>
&=
\frac{1}{2\pi} \int_\RR  \frac{c_o}{\omega} \hat{F}(\omega) {\rm Im} \hat{G} \big(\omega,\bx,\bx') \big) 
 \exp\big[ i \omega  ( t-t' ) \big]
 d\omega   \\
 &=
\frac{1}{8\pi^2} \int_\RR  \hat{F}(\omega) {\rm sinc} \Big(\frac{\omega}{c_o} \big|\bx- \bx'\big| \Big) 
 \exp\big[ i \omega  ( t-t' ) \big]
 d\omega ,
  \label{eq:exprescov3do}
\end{align}
where ${\rm sinc}(x)=\sin(x)/x$.
This result is in fact very general and holds true for a large class of spatial source distributions,
as was shown in the literature about seismic interferometry \cite{curtis06} or ambient noise imaging \cite{garpapa16}.
{From now on we will assume that the spatial distribution of the noise sources is such that (\ref{eq:exprescov3do}) holds true}.

The random process $(u(t,\bx))_{t \in \RR,\bx\in \RR^3} $ has Gaussian statistics, mean zero, and covariance function (\ref{eq:exprescov3do}).
An equivalent description of such a random field is:
\begin{equation}
\label{eq:represboyer}
u(t,\bx) = 
\int_{\RR^3} a(\bk) \exp\big[ i (\bk\cdot\bx -  \Omega(\bk) t)\big] d\bk ,
\end{equation}
where $\Omega(\bk)=c_o |\bk|$ is the dispersion relation of the three-dimensional
wave equation and $(a(\bk))_{\bk \in \RR^3}$ is a complex-valued Gaussian process with mean zero
and covariance function
\begin{equation}
\left< a(\bk) \overline{a(\bk')}\right> = {\cal A}(|\bk|) \delta(\bk-\bk') ,
\end{equation}
with $a(-\bk) = \overline{a(\bk)}$ and 
\begin{equation}
\label{eq:calAk}
{\cal A}(k) = \frac{c_o}{32\pi^3 k^2} \hat{F}(c_ok).
\end{equation}
The proof of the equivalence is straightforward: the random process defined by (\ref{eq:represboyer})
has Gaussian statistics (as it is a linear transform of a Gaussian process), mean zero, and its covariance
is
\begin{align*}
\left< u(t,\bx) u(t',\bx') \right>
&=
\int_{\RR^3} \exp\big[ i \big( \bk \cdot (\bx-\bx') -  \Omega(\bk) (t-t')\big)\big]  {\cal A}(|\bk|) d\bk\\
&=
\int_0^\infty  \int_{\mathbb{S}^3}  \exp\big[ i \big( k \hat{\bk} \cdot (\bx-\bx') \big]  d \hat{\bk} \,
{\cal A}(k)  \exp \big[ -i  c_o k (t-t')\big] k^2 dk \\
&= 4\pi 
\int_0^\infty {\rm sinc}( k|\bx-\bx'|)
k^2 {\cal A}(k)  \exp \big[- i c_o k (t-t')\big] dk ,
\end{align*}
that is equal to (\ref{eq:exprescov3do}) after the change of variable $k=\omega/c_o$.
This shows that $u(t,\bx)$ can be considered a superposition of uncorrelated plane waves. 
We are in a situation similar to the one addressed in \cite{boyer80} and we will recover the results of this paper when the medium is homogeneous.
However, our approach makes it possible to address non-homogeneous media as we will see in Section \ref{sec:2}.

\begin{figure}
\begin{center}
\begin{tabular}{c}
\includegraphics[width=6.4cm]{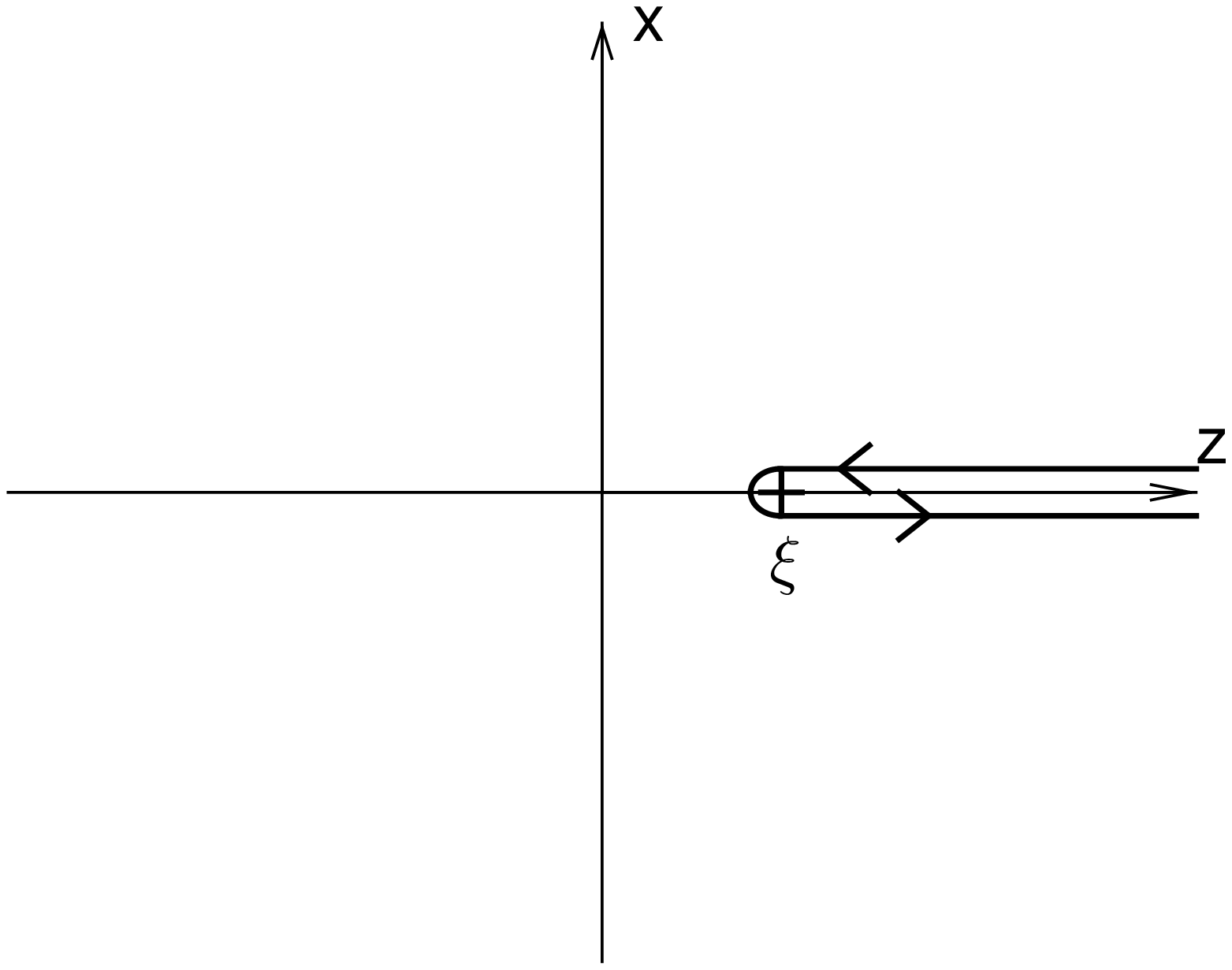}
\includegraphics[width=6.4cm]{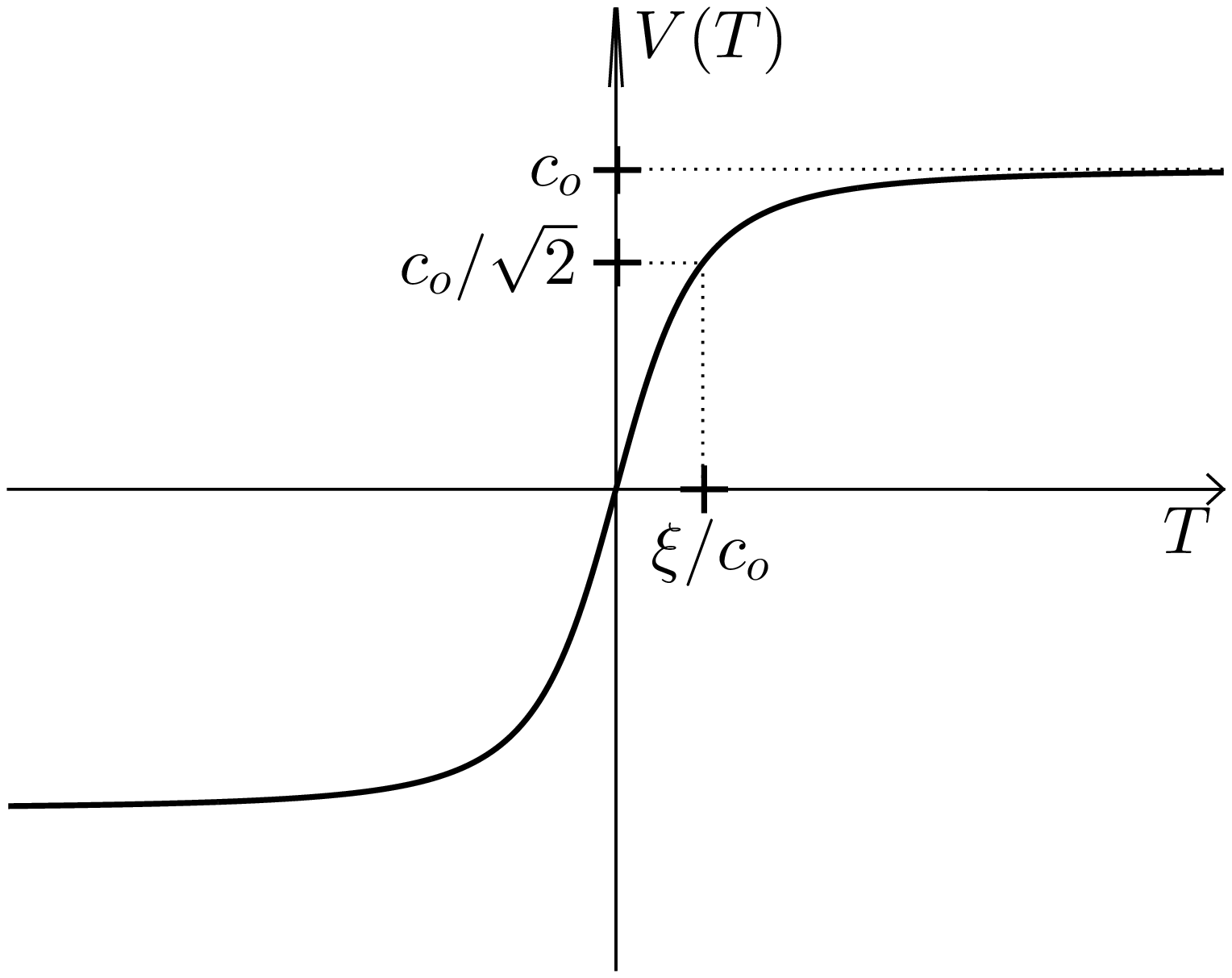}
\end{tabular}
\end{center}
\caption{
Rindler trajectory (\ref{eq:rindler3d}) in a three-dimensional open medium (left) and  speed (along the $z$-axis) of the Rindler observer in the laboratory frame (right).
    \label{fig:rindler1}
}
\end{figure}

\subsection{The Wigner transform of the recorded signal}

We consider an observer with a Rindler trajectory (see Figure \ref{fig:rindler1}, left)
whose time-space coordinates are (with fixed $\xi>0$):
\begin{equation}
\label{eq:rindler3d}
T(\tau) = ({\xi}/{c_o}) \sinh( c_o \tau / \xi) ,\quad \quad \bX(\tau) = \big(0,0,\xi\cosh( c_o \tau / \xi)\big)  ,
\end{equation}
where $\tau$ is the proper time (the  time as perceived by the accelerated observer).
The Rindler trajectory (\ref{eq:rindler3d}) describes a trajectory with constant acceleration \cite{rindler,book:ulf},
in the sense that the acceleration of the observer relative to its instantaneous inertial rest frame is constant and equal to $g=c_o^2/\xi$.
In the laboratory frame the trajectory is along the $z$-axis with the coordinate $Z(T)=\sqrt{\xi^2+c_o^2 T^2}$ and  the speed $V(T)=\partial_T Z(T)=
c_o^2 T / \sqrt{\xi^2+c_o^2 T^2}$  (see Figure \ref{fig:rindler1}, right).
We will see in the next subsection that this trajectory satisfies a unique property from the point of view of correlation of noise. 

The signal recorded by the Rindler observer is
\begin{equation}
U(\tau) = u \big( T( \tau ),  \bX(\tau)\big)  .
\end{equation}
We look for the properties of the autocorrelation function of the recorded signal:
\begin{equation}
\label{def:auto}
\left< U(\tau+\frac{\tau'}{2})U(\tau-\frac{\tau'}{2})\right> ,
\end{equation}
and the Wigner transform of the recorded signal, that is, its local spectrum:
\begin{equation}
\label{def:wigner}
W(\tau,\omega)= \int_\RR \left< U(\tau+\frac{\tau'}{2})U(\tau-\frac{\tau'}{2})\right>  \exp(i \omega \tau' )d\tau'  .
\end{equation}
From (\ref{eq:exprescov3do}) the autocorrelation function of the noise recorded by the Rindler observer is
\begin{align}
\nonumber
&\left< U(\tau+\frac{\tau'}{2})U(\tau-\frac{\tau'}{2})\right>\\
\nonumber
&=
\frac{1}{2\pi} \int_\RR  \frac{c_o}{\omega} \hat{F}(\omega) {\rm Im} \hat{G} \Big(\omega,\bX(\tau+\frac{\tau'}{2}),\bX(\tau-\frac{\tau'}{2}) \Big) 
 \exp\Big[ i \omega \big( T(\tau+\frac{\tau'}{2})-T(\tau-\frac{\tau'}{2})\big) \Big]
 d\omega   \\
 &=
\frac{1}{8\pi^2} \int_\RR  \hat{F}(\omega) {\rm sinc} \Big(\frac{\omega}{c_o} \big|\bX(\tau+\frac{\tau'}{2})- \bX(\tau-\frac{\tau'}{2})\big| \Big) 
 \exp\Big[ i \omega \big( T(\tau+\frac{\tau'}{2})-T(\tau-\frac{\tau'}{2})\big) \Big]
 d\omega .
 \label{eq:expresauto3do}
\end{align}

We will focus our attention to the case where the source spectrum is of the form
\begin{equation}
\label{eq:formsourcespectrum}
\hat{F}(\omega)= f_o |\omega| .
\end{equation}
We explain in the next subsection the special properties of this source spectrum.
{Note that, in practice, the source spectrum can be a finite-energy approximation of the ideal
spectrum (\ref{eq:formsourcespectrum}), such as $\hat{F}(\omega)=f_o |\omega|\exp(-\epsilon |\omega|)$ for some $\epsilon>0$ for instance. The following results are then valid in the sense explained in Remark \ref{rem:1}.}

If the source spectrum is (\ref{eq:formsourcespectrum}) then the power spectral density of the signal 
recorded by a stationary observer is $f_o |\omega| / (4\pi)$ by (\ref{eq:exprescov3do}) (with $\bx'=\bx$).
For an observer with a Rindler trajectory,
we find from (\ref{eq:expresauto3do}) that, for any $\tau' \neq 0$
 {(see \ref{app:B1})}:
\begin{align}
\left< U(\tau+\frac{\tau'}{2})U(\tau-\frac{\tau'}{2})\right>
=
- \frac{c_o^2 f_o}{16\pi^2 \xi^2} \frac{1}{\sinh^2(c_o \tau'/(2\xi))},
\label{eq:expresauto3dot}
\end{align} 
so that the Wigner transform is independent of $\tau$
 {(see \ref{app:B2})}:
\begin{align}
W(\tau,\omega) = \frac{f_o}{4\pi} \frac{\omega}{\tanh(\pi \xi \omega/c_o)}.
\label{def:Wonu0}
\end{align}
At any time $\tau$, the observer feels the same spectrum.
Moreover, the spectrum is a perturbation of the spectrum $f_o|\omega|/(4\pi)$ observed by a stationary observer and it has the ``Planck" form:
\begin{equation}
W(\tau,\omega) = W_o(\omega), \quad \quad W_o(\omega) =
\frac{f_o |\omega|}{4\pi} \Big( 1+\frac{2}{e^{2\pi \xi |\omega|/c_o}-1} \Big).
\label{def:Wonu}
\end{equation}
The analogy with the Planck spectrum is obtained by identifying $2\pi \xi / c_o$ and $\bar{h}/(K_BT)$ 
(with $\bar{h}$ the Planck's constant divided by $2\pi$, $K_B$ the Boltzmann's constant, and $T$ the Unruh temperature \cite{unruh}).
The result (\ref{def:Wonu}) was obtained from the representation (\ref{eq:represboyer}) in \cite{boyer80,boyer84}.
This representation was convenient to analyze the system in an open medium,
but the presence of an obstacle imposes the formulation with noise sources.

{
\begin{remark}
\label{rem:1}
In practice the spectrum may not be equal to (\ref{eq:formsourcespectrum}) but of the form
$\hat{F}(\omega)=f_o |\omega|\exp(-\epsilon |\omega|)$ for some $\epsilon>0$ for instance.
This situation is analyzed in \ref{app:B3}.
The autocorrelation function  of the recorded signal is then given 
by (\ref{eq:Ueps}) for any $\epsilon >0$ which reduces to (\ref{eq:expresauto3dot}) when $\epsilon$ is small and $|\tau'|$ is larger than $O(\epsilon)$.
The  Wigner transform of the recorded signal is given  by (\ref{eq:Weps}) for any $\epsilon >0$
which reduces to (\ref{def:Wonu0}) when $\epsilon$ is small and $\omega$ is smaller than $O(\epsilon^{-1})$.
This remark shows that our results are somewhat robust with respect to the form of the source spectrum
and that it is possible to consider a source spectrum with finite energy and amplitude while preserving the main results.
\end{remark}
}

{
\begin{remark}
\label{rem:2}
In order to compute the Wigner transform (\ref{def:wigner}) it is necessary to integrate over $\RR$.
If the signal is only recorded over a finite time interval, then one can compute 
$$
W_\chi(\tau,\omega) = \int_\RR \chi(\tau') \left< U(\tau+\frac{\tau'}{2})U(\tau-\frac{\tau'}{2})\right>  \exp(i \omega \tau' )d\tau'  ,
$$
where $\chi$ is a cut-off function. We then have
$$
W_\chi(\tau,\omega) = \int W(\tau,\omega-\omega') \hat{\chi}(\omega') d\omega' ,
$$
where $\hat{\chi}$ is the Fourier transform of $\chi$.
If the recording time interval has duration $T_c$, then this means that we can extract the Wigner transform up to a convolution
in $\omega$ with a kernel of width $1/T_c$.
This remark shows that our results are somewhat robust with respect to the duration of the recording time interval 
and that it is possible to consider a finite recording time interval  while preserving the main results.
\end{remark}
}

\subsection{Properties of the source and signal spectra}
{We recall that we assume that the spatial distribution of the noise sources is such that (\ref{eq:exprescov3do}) holds true.}
The following proposition justifies why we focus our attention to the case where the source spectrum has the form (\ref{eq:formsourcespectrum}).
\begin{proposition}
The source spectrum $\hat{F}(\omega)=|\omega|$ (up to a multiplicative constant) is the only Lorentz-invariant spectrum.
\end{proposition}

This was  proved in \cite{boyer80}.
If the source spectrum is Lorentz-invariant, then the signal recorded by an observer with constant velocity has the same spectrum,
whatever the velocity of the observer (provided it is constant). 
This result can be recovered by using our formulation of the problem, except that two possible forms of source spectrum are possible, {as stated in the following proposition proved in \ref{proof:prop2}.}

\begin{proposition}
\label{prop:2}
The only source spectra for which  the signal recorded by an observer with constant velocity has a spectrum independent of  velocity
are of the form $\hat{F}(\omega)=f_o |\omega| +f_1 /|\omega|$.
\end{proposition}

The following proposition {proved in  \ref{proof:prop3a}}
underlines an important property of a source spectrum of the form (\ref{eq:formsourcespectrum}).
\begin{proposition}
\label{prop:3a}
The source spectrum $\hat{F}(\omega)=|\omega|$ (up to a multiplicative constant) is the only one
for which the recorded signal is stationary (i.e., its autocorrelation function does not depend on $\tau$,
or equivalently its Wigner transform does not depend on $\tau$) for a Rindler trajectory.
\end{proposition}

The following proposition 
clarifies a unique and crucial property of the Rindler trajectory.

\begin{proposition}
\label{prop:3}%
The Rindler trajectories are the only straight-line trajectories 
for which the recorded signal is stationary (i.e., its autocorrelation function does not depend on $\tau$,
or equivalently its Wigner transform does not depend on $\tau$) when the power spectral density
of the noise sources is Lorentz invariant (i.e. $\hat{F}(\omega)= f_o |\omega|$).
\end{proposition}

{\it Proof.}
Let us consider a trajectory with proper time $\tau$ of the form
$(t(\tau),{\bx}(\tau))$ with ${\bx}(\tau)=(0,0,z(\tau))$. The time in the laboratory frame $t$ is related to 
the proper time $\tau$  by (\ref{eq:mink}).
By (\ref{eq:exprescov3do}) 
the autocorrelation function of the recorded signal has the form
\begin{align}
\nonumber
&\left< U(\tau+\frac{\tau'}{2})U(\tau-\frac{\tau'}{2})\right>
\\
\nonumber
 &=
\frac{f_o}{4\pi^2} \int_0^\infty
\frac{{\rm sin} \Big(\frac{\omega}{c_o} \big| z(\tau+\frac{\tau'}{2})- z(\tau-\frac{\tau'}{2})\big| \Big) }{
\frac{1}{c_o} \big| z(\tau+\frac{\tau'}{2})-  z(\tau-\frac{\tau'}{2})\big| }
 \cos\Big( \omega \big( t(\tau+\frac{\tau'}{2})-t(\tau-\frac{\tau'}{2})\big) \Big)
 d\omega \\
 &=
 \frac{f_o}{4\pi^2} \frac{1}{\frac{1}{c_o^2} \big( z(\tau+\frac{\tau'}{2})-  z(\tau-\frac{\tau'}{2})\big)^2 - 
 \big( t(\tau+\frac{\tau'}{2})-t(\tau-\frac{\tau'}{2})\big)^2}.
\end{align}
On the one hand,
if the trajectory is Rindler, then there exist $t_0,z_0,\tau_0,\xi$ such that $t(\tau)=t_0+ ({\xi}/{c_o}) \sinh( c_o (\tau-\tau_0) / \xi)$ and $z(\tau) = z_0+\xi\cosh( c_o (\tau-\tau_0) / \xi)\big) $,
so that 
\begin{align*}
t(\tau+\frac{\tau'}{2})-t(\tau-\frac{\tau'}{2}) =&2 ({\xi}/{c_o}) \cosh (c_o (\tau-\tau_0) / \xi) \sinh(c_o \tau'/ (2c_o)),\\
z(\tau+\frac{\tau'}{2})-  z(\tau-\frac{\tau'}{2}) =&2 {\xi}  \sinh (c_o (\tau-\tau_0) / \xi) \sinh(c_o \tau'/ (2c_o)),
\end{align*}
and therefore
$$
\frac{1}{c_o^2} \big( z(\tau+\frac{\tau'}{2})-  z(\tau-\frac{\tau'}{2})\big)^2 - 
 \big( t(\tau+\frac{\tau'}{2})-t(\tau-\frac{\tau'}{2})\big)^2
 = - \frac{4\xi^2}{c_o^2} \sinh^2(c_o \tau'/ (2c_o)),
$$
which indeed implies  that the autocorrelation function of the recorded signal is independent of $\tau$.

On the other hand,
if we impose that the autocorrelation function of the recorded signal is independent of $\tau$, then this means that 
$$
\big( z(\tau+\tau')-z(\tau) \big)^2 - \Big( \int_{\tau}^{\tau+\tau'} \sqrt{c_o^2 + |\dot{z}(\tau'')|^2} d\tau'' \Big)^2
$$
should be independent of $\tau$, where we have used (\ref{eq:mink}).
We can expand this expression for small $\tau'$ and we get that the fourth-order coefficient of the Taylor series expansion
is equal to
$$
-\frac{1}{12} \frac{c_o^2 \ddot{z}^2}{c_o^2+\dot{z}^2} .
$$
This coefficient should be independent of $\tau$, which means that the function ${V}=\dot{z}/c_o$ should satisfy
an ordinary differential equation of the form
$$
\frac{\dot{{V}}}{\sqrt{1+{V}^2}} =c_1 ,
$$
for some constant $c_1$.
The general solution of this equation is
$$
{V}(\tau) = \sinh (c_1\tau+c_2),
$$
which corresponds to $z(\tau) = c_o/c_1 \cosh(c_1\tau+c_2)+ c_3$ and $t(\tau) = (1/c_1) \sinh(c_1\tau+c_2)+ c_4$,
sor some constants $c_1,\ldots,c_4$.
We conclude that the trajectory should be a Rindler trajectory.
\qed

Finally the following proposition shows that the straight-line trajectory hypothesis in Proposition \ref{prop:3} is important.
\begin{proposition}
\label{prop:4}
There exist trajectories that do not follow a straight line but that  give stationary recorded signals when the power spectral density
of the noise sources is of the form $\hat{F}(\omega)= f_o |\omega|$.
\end{proposition}

{\it Proof.}
Let us consider a uniform circular motion with proper time $\tau$ of the form
$(t(\tau),{\bx}(\tau))$ with ${\bx}(\tau)=(x(\tau),y(\tau),0)$,
$$
t=\gamma \tau, 
\quad x(\tau)=\frac{c_o\sqrt{\gamma^2-1}}{p} \cos(p\tau),
\quad y(\tau)=\frac{c_o\sqrt{\gamma^2-1}}{p} \sin(p\tau).
$$
Here $\gamma>1$ is the Lorentz factor, $p=\gamma p_0$, $p_0$ is the coordinate angular velocity.
The time in the laboratory frame $t$ is related to 
the proper time $\tau$ by
$\dot{t}^2 - (\dot{x}^2 +\dot{y}^2)/c_o^2=1$.
Then we find that 
\begin{align}
\left< U(\tau+\frac{\tau'}{2})U(\tau-\frac{\tau'}{2})\right> 
= \frac{f_o}{4\pi^2} \frac{1}
 {4 \frac{\gamma^2-1}{p^2} \sin^2\big(\frac{p\tau'}{2} \big)
 -\gamma^2 {\tau'}^2}
 ,
 \end{align}
which does not depend on $\tau$.
The spectrum $W(\tau,\omega)$ is independent of $\tau$ as well, it is a perturbation of the spectrum observed by a stationary observer:
$$
W(\tau,\omega) = \frac{f_o |\omega|}{4\pi} +\frac{f_o p}{4\pi} W_\gamma\big( \frac{\omega}{p}\big) ,
$$
with
$$
W_\gamma(w) = \frac{\gamma^2-1}{4\pi^2}  \int_\RR
\frac{s^2-\sin^2(s) }{s^2 (\gamma^2s^2-(\gamma^2-1) \sin^2(s) )}  \cos(2 w s) ds ,
$$
which is an even, bounded, and integrable function.
The perturbation {does not have a Planck form}.
For $0\leq \gamma^2-1 \ll 1$, we have
\begin{align*}
W_\gamma(w) &=
 \frac{\gamma^2-1}{4\pi^2}  \int_\RR
\frac{s^2-\sin^2(s) }{s^4}  \cos(2 w s) ds
+o(\gamma^2-1) \\
& =
 \frac{\gamma^2-1}{6\pi} (1-|w|)_+^3 +o(\gamma^2-1).
\end{align*}
\qed

It would be interesting to identify all trajectories that give stationary recorded signals,
but this is beyond the scope of this paper.
We may conjecture that the result should be that the acceleration should be constant, as it is for the Rindler trajectory and the circular trajectory,
but it is not so straightforward.
Indeed, if the trajectory ${\bx}(\tau)=(x(\tau),y(\tau),z(\tau))$ gives a stationary recorded signal, then
the autocorrelation function of the recorded signal  
\begin{align}
\nonumber
\left< U(\tau+\frac{\tau'}{2})U(\tau-\frac{\tau'}{2})\right>
= \frac{f_o}{4\pi^2} \frac{1}{\frac{1}{c_o^2} \big|\bx(\tau+\frac{\tau'}{2})-\bx(\tau-\frac{\tau'}{2})\big|^2  - 
 \big( t(\tau+\frac{\tau'}{2})-t(\tau-\frac{\tau'}{2})\big)^2}
\end{align}
should be independent of $\tau$,
in other words, 
$$
\big|\bx(\tau+\tau')-\bx(\tau)\big|^2 - \Big( \int_{\tau}^{\tau+\tau'} \sqrt{c_o^2 + |\dot{\bx}(\tau'')|^2} d\tau'' \Big)^2
$$
should be independent of $\tau$. By following the same strategy as in the proof of Proposition \ref{prop:3},
we find that the normalized velocity ${\itbf V}=\dot{\bx}/c_o$ should satisfy a nonlinear ordinary differential equation of the form 
$$
|\dot{{\itbf V}}|^2 - \frac{({\itbf V}\cdot\dot{{\itbf V}})^2}{1+|{\itbf V}|^2} =c_1,
$$
for some constant $c_1$.
However, it does not seem straightforward to identify all the solutions of this equation.
Rindler trajectories and uniform circular motions are solutions, but there are other solutions, such as constant helicoidal motions:
$$
(x(\tau),y(\tau)) = \frac{c_o \sqrt{\gamma^2-1}\sqrt{\alpha}}{p} \big(
\cos (p\tau),\sin (p\tau) \big) ,\quad
z(\tau) = c_o \sqrt{\gamma^2-1} \sqrt{1-\alpha} \tau, 
$$
and $t(\tau)=\gamma\tau$, where $\tau$ is the proper time and $\gamma>1$, $p>0$, $\alpha \in [0,1]$ are constant parameters,
or uniformly accelerated helicoidal motions:
$$
(x(\tau),y(\tau)) = \frac{c_o A}{p} \big(
\cos (p\tau),\sin (p\tau) \big) ,\quad
z(\tau) = \xi \sqrt{A^2+1} \cosh(c_o \tau /\xi),
$$
and $t(\tau)=\sqrt{A^2+1}(\xi/c_o)  \sinh(c_o \tau/\xi)$, where $\tau$ is the proper time and $A ,\xi,p$  are constant parameters
{(we say that this motion is uniformly accelerated because the acceleration along the $z$-axis of the observer relative to its instantaneous inertial rest frame is constant and equal to $(1+A^2) c_o^2 / \xi$, but the acceleration in the $x$- and $y$-directions are not constant)}.
These trajectories also satisfy the property that the signal recorded by a moving observer is stationary.

\section{A Rindler observer in front of an obstacle}
\label{sec:2}
\subsection{The configuration}
We consider the three-dimensional scalar wave equation for the scalar wave field $u(t,\bx)$ in the half-space $\RR^2 \times(0,+\infty)$:
\begin{equation}
\frac{1}{c_o^2} \partial_t^2 u -\Delta u =n(t,\bx) ,
\end{equation}
with a reflecting (Dirichlet) boundary condition on the plane $z=0$ (with $\bx=(x,y,z)$):
\begin{equation}
u(t,(x,y,0))=0 .
\end{equation}
The goal of this section is to show that an observer following a Rindler trajectory can detect and localize
the obstacle (i.e. the interface $z=0$) from the signal that it records.

The source term $n(t,\bx)$ models a noise source distribution. It is a zero-mean process, stationary in time
and delta-correlated in space:
\begin{equation}
\left< n(t,\bx)n(t',\bx') \right> = F(t-t')\delta(\bx-\bx') K(\bx) ,
\end{equation}
where $K(\bx)$ is the function that characterizes the spatial support of the noise source distribution
(it is supported in the half-space $\RR^2 \times(0,+\infty)$)
and  the Fourier transform $\hat{F}(\omega)$ of $F$ is the power spectral density of the sources.
For simplicity (this can be easily generalized)
we can consider that the function $K$ is spatially supported on a surface of a half-ball $\partial B_+({\bf 0},L) = 
\{ \bx =(x,y,z)\in \RR^3 \, , \, z>0, \, |\bx|=L\}$  with a very large radius $L$:
$$
K(\bx) = \delta_{\partial B_+({\bf 0},L)} (\bx) .
$$

\begin{figure}
\begin{center}
\begin{tabular}{c}
\includegraphics[width=6.4cm]{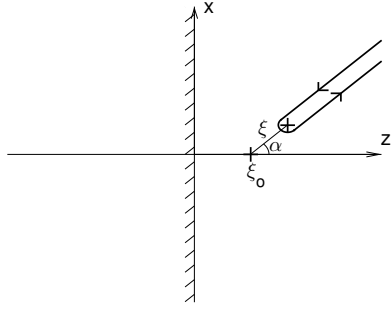}
\end{tabular}
\end{center}
\caption{Rindler trajectory (\ref{eq:rindlerobs}) with an obstacle in the plane $z=0$.
    \label{fig:rindler2}
}
\end{figure}

We consider an observer with a Rindler trajectory whose time-space coordinates are
 (with fixed $\xi>0$,  $\alpha \in (-\pi/2,\pi/2)$, $\xi_o> - \xi\cos\alpha$):
\begin{equation}
\label{eq:rindlerobs}
T(\tau) = ({\xi}/{c_o}) \sinh(c_o\tau/\xi) ,\quad \quad \bX(\tau) = (0,0,\xi_o)+\xi \cosh(c_o\tau/\xi) \big(\sin \alpha , 0, \cos \alpha\big) ,
\end{equation}
where $\tau$ is the proper time of the observer. 
With this parameterization, the position of the observer is the closest to the obstacle at $\tau=0$, 
when it is located at $\bX(0)=(\xi \sin \alpha , 0 ,\xi_o +\xi  \cos \alpha)$ (see Figure \ref{fig:rindler2}). 
We will first address the direct problem: calculation of the autocorrelation function of the field recorded by the observer.
In the inverse problem, the observer knows its 
acceleration $c_o^2/\xi$ 
and its proper time $\tau$, 
it can observe the  autocorrelation function of the recorded field and it looks for its relative position to the obstacle, that is to say, it looks for $\alpha$ and $\xi_o$.

\subsection{The Wigner transform of the recorded signal}
We first derive an expression of the autocorrelation function in the presence of the obstacle (the interface $z=0$).

\begin{proposition}
The autocorrelation function of the recorded signal has the form
\begin{align}
\nonumber
&
\left< U(\tau+\frac{\tau'}{2})U(\tau-\frac{\tau'}{2})\right>
\\
\nonumber
&=
\frac{1}{8\pi^2} \int_\RR  \hat{F}(\omega) {\rm sinc} \Big(\frac{\omega}{c_o} \big|\bX(\tau+\frac{\tau'}{2})- \bX(\tau-\frac{\tau'}{2})\big| \Big) 
 \exp\Big[ i \omega \big( T(\tau+\frac{\tau'}{2})-T(\tau-\frac{\tau'}{2})\big) \Big]
 d\omega \\
&\quad - \frac{1}{8\pi^2} \int_\RR  \hat{F}(\omega) {\rm sinc} \Big(\frac{\omega}{c_o}\big|\bX^{\rm s}(\tau+\frac{\tau'}{2})- \bX(\tau-\frac{\tau'}{2})\big| \Big) 
 \exp\Big[ i \omega \big( T(\tau+\frac{\tau'}{2})-T(\tau-\frac{\tau'}{2})\big) \Big]
 d\omega  ,
  \label{eq:expresauto3d2}
\end{align}
with 
$\bX^{\rm s}(\tau)=  (0,0,-\xi_o)+\xi \cosh(c_o \tau/\xi) \big(0,\sin \alpha , -\cos \alpha\big)$.
\end{proposition}

{\it Proof.}
The expression of the autocorrelation function is
\begin{align*}
\left< U(\tau+\frac{\tau'}{2})U(\tau-\frac{\tau'}{2})\right>
=&
\frac{1}{2\pi} \int_\RR \int_{\RR^3} \overline{\hat{\cal G}\Big(\omega,\bX(\tau+\frac{\tau'}{2}),\by\Big)}
\hat{\cal G}\Big(\omega,\bX(\tau-\frac{\tau'}{2}),\by\Big) \\
&\quad \times K(\by) \hat{F}(\omega) \exp\Big[ i \omega \big( T(\tau+\frac{\tau'}{2})-T(\tau-\frac{\tau'}{2})\big) \Big]
d\by d\omega  .
\end{align*}
Here $\hat{\cal G}(\omega,\bx,\by)$ is the Green's function in the presence of the reflecting plane,
that is to say
$$
\hat{\cal G}(\omega,\bx,\by) = \hat{G}(\omega,\bx,\by)- \hat{G}(\omega,\bx, \by^{\rm s} ) ,
$$
where $\by^{\rm s}= (y_1,y_2,-y_3)$ is the symmetric point of $ \by = (y_1,y_2,y_3)$ and $\hat{G}(\omega,\bx,\by)$
is the three-dimensional homogeneous Green's function (\ref{eq:green3d}).
Using the fact that $\hat{G}(\omega,\bx,\by^{\rm s})=\hat{G}(\omega,\bx^{\rm s},\by)$,
we find that 
\begin{align*}
\left< U(\tau+\frac{\tau'}{2})U(\tau-\frac{\tau'}{2})\right>
=&
\frac{1}{2\pi}  \int_\RR \int_{\RR^3} \overline{\hat{G}\Big(\omega,\bX(\tau+\frac{\tau'}{2}),\by\Big)}
\hat{G}\Big(\omega,\bX(\tau-\frac{\tau'}{2}),\by\Big) \\
&\quad \times  \hat{F}(\omega) K^{\rm s}(\by) \exp\Big[ i \omega \big( T(\tau+\frac{\tau'}{2})-T(\tau-\frac{\tau'}{2})\big) \Big]
d\by d\omega \\
& - \frac{1}{2\pi}  \int_\RR \int_{\RR^3} \overline{\hat{G}\Big(\omega,\bX^{\rm s}(\tau+\frac{\tau'}{2}),\by\Big)}
\hat{G}\Big(\omega,\bX(\tau-\frac{\tau'}{2}),\by\Big) \\
&\quad \times  \hat{F}(\omega) K^{\rm s}(\by) \exp\Big[ i \omega \big( T(\tau+\frac{\tau'}{2})-T(\tau-\frac{\tau'}{2})\big) \Big]
d\by d\omega ,
\end{align*}
where we have defined
$$
K^{\rm s}(\by) = K(\by)+K(\by^{\rm s})  .
$$
As we assume that  $K(\bx) = \delta_{\partial B_+({\bf 0},L)} (\bx)$,
 $K^{\rm s}$ is supported at the surface of the ball with center at ${\bf 0}$ and radius $L$, i.e. $K^{\rm s}(\bx) = \delta_{\partial B({\bf 0},L)} (\bx)$,
so  we can invoke Helmholtz-Kirchhoff identity~(\ref{hk})
to compute the integral in $\by$ and we get
\begin{align}
\nonumber
&
\left< U(\tau+\frac{\tau'}{2})U(\tau-\frac{\tau'}{2})\right>
\\
\nonumber
&=
\frac{1}{2\pi} \int_\RR  \frac{c_o}{\omega} \hat{F}(\omega) {\rm Im} \hat{G} \Big(\omega,\bX(\tau+\frac{\tau'}{2}),\bX(\tau-\frac{\tau'}{2}) \Big) 
 \exp\Big[ i \omega \big( T(\tau+\frac{\tau'}{2})-T(\tau-\frac{\tau'}{2})\big) \Big]
 d\omega \\
&\quad - \frac{1}{2\pi} \int_\RR  \frac{c_o}{\omega}  \hat{F}(\omega) 
{\rm Im} \hat{G}\Big(\omega,\bX^{\rm s}(\tau+\frac{\tau'}{2}),\bX(\tau-\frac{\tau'}{2}) \Big) 
\exp\Big[ i \omega \big( T(\tau+\frac{\tau'}{2})-T(\tau-\frac{\tau'}{2})\big) \Big]
 d\omega  .
 \label{eq:cross2}
\end{align}
We can rewrite this equation in two explicit forms, either (\ref{eq:expresauto3d2}) or 
\begin{align}
\nonumber
&
\left< U(\tau+\frac{\tau'}{2})U(\tau-\frac{\tau'}{2})\right>
\\
\nonumber
&=
\frac{1}{8\pi} \int_{-1}^1 F \Big(\frac{v}{c_o} \big|\bX(\tau+\frac{\tau'}{2})- \bX(\tau-\frac{\tau'}{2})\big| +
\big( T(\tau+\frac{\tau'}{2})-T(\tau-\frac{\tau'}{2})\big) \Big)
 dv \\
&\quad - \frac{1}{8\pi} \int_{-1}^1 F \Big( \frac{v}{c_o} \big|\bX^{\rm s}(\tau+\frac{\tau'}{2})- \bX(\tau-\frac{\tau'}{2}) \big|+
 \big( T(\tau+\frac{\tau'}{2})-T(\tau-\frac{\tau'}{2}\big) \Big)
 dv  .
   \label{eq:expresauto3d2b}
\end{align}
{The form (\ref{eq:expresauto3d2b}) is obtained by using the explicit form (\ref{eq:green3d}) of the Green's function and the identity $ \int_{-1}^1\exp(i v s) dv=2{\rm sinc}( s)$, so that 
\begin{align*}
&\int_\RR  \frac{c_o}{\omega} \hat{F}(\omega) {\rm Im} \hat{G} \Big(\omega,\bX(\tau+\frac{\tau'}{2}),\bX(\tau-\frac{\tau'}{2}) \Big) 
 \exp\Big[ i \omega \big( T(\tau+\frac{\tau'}{2})-T(\tau-\frac{\tau'}{2})\big) \Big]
 d\omega \\
& =  \frac{1}{4\pi} \int_\RR \hat{F}(\omega) 
{\rm sinc}  \Big( \frac{\omega}{c_o} \big|\bX(\tau+\frac{\tau'}{2}) - \bX(\tau-\frac{\tau'}{2}) \big| \Big) 
\exp \Big[ i   \omega \big( T(\tau+\frac{\tau'}{2})-T(\tau-\frac{\tau'}{2})\big) \Big]d\omega 
\\
& = \frac{1}{8\pi} \int_{-1}^1  dv\int_\RR \hat{F}(\omega) 
 \exp \Big[ i   \frac{\omega v}{c_o} \big|\bX(\tau+\frac{\tau'}{2}) - \bX(\tau-\frac{\tau'}{2}) \big| \Big]
\exp \Big[ i   \omega \big( T(\tau+\frac{\tau'}{2})-T(\tau-\frac{\tau'}{2})\big) \Big]d\omega 
\\
&=\frac{1}{4}
 \int_{-1}^1 dv F \Big(\frac{v}{c_o} \big|\bX(\tau+\frac{\tau'}{2})- \bX(\tau-\frac{\tau'}{2})\big| +
\big( T(\tau+\frac{\tau'}{2})-T(\tau-\frac{\tau'}{2})\big) \Big) ,
\end{align*}
and similarly for the second term in (\ref{eq:cross2}).
}
\qed

From now on we assume that the power spectral density of the sources is of the form (\ref{eq:formsourcespectrum}).
If the observer is stationary at position $\bX_0$, then we have from (\ref{eq:expresauto3d2}) taken with $\bX(\tau)\equiv \bX_0$ and $T(\tau)\equiv \tau$:
\begin{align*}
\left< u(\tau +\frac{\tau'}{2},\bX_0) u(\tau -\frac{\tau'}{2},\bX_0) \right>
&=
\frac{1}{8\pi^2} \int_\RR f_o |\omega|\exp ( i \omega \tau')
 d\omega \\
&\quad - \frac{1}{8\pi^2} \int_\RR  f_o |\omega| {\rm sinc} \Big(\frac{\omega}{c_o}\big|\bX^{\rm s}_0 -  \bX_0\big| \Big) 
 \exp( i \omega\tau')
 d\omega  ,
\end{align*}
and therefore
\begin{align*}
 \int_\RR \left< u(\tau +\frac{\tau'}{2},\bX_0) u(\tau -\frac{\tau'}{2},\bX_0) \right> \exp(i \omega \tau' )d\tau'  =
 \frac{f_o |\omega|}{4\pi}\Big( 1 -  {\rm sinc} \Big(\frac{\omega}{c_o}\big|\bX^{\rm s}_0 -  \bX_0\big| \Big)  \Big) .
\end{align*}
This shows that the unperturbed spectrum $f_o |\omega|/(4\pi)$ felt by a stationary observer in a homogeneous medium
is perturbed by the obstacle and that the observer
can extract the distance $|\bX^{\rm s}_0 -  \bX_0|$ from the spectrum,
that is to say, twice the distance from the observer to the obstacle.
However, the observer cannot determine the angular position of the obstacle, which is not surprising by symmetry of the system.

Let us now consider an observer on a Rindler trajectory.
By (\ref{eq:expresauto3d2}) we have 
\begin{align}
\nonumber
\left< U(\tau+\frac{\tau'}{2})U(\tau-\frac{\tau'}{2})\right>
=
&- \frac{c_o^2 f_o}{16\pi^2 \xi^2}\frac{1}{\sinh^2(\eta'/2)} 
\\
&+ \frac{c_o^2 f_o}{16 \pi^2 \xi^2}
\frac{1}{A \cosh^2(\eta'/2) +B(\eta)  \cosh(\eta'/2)+C(\eta)}  ,
\label{eq:automur}
\end{align}
where $\eta= c_o \tau /\xi$ and $\eta'=c_o\tau'/\xi$,
\begin{equation}
\label{eq:defABC}
A= \sin^2\alpha,\quad \quad B(\eta)= -2\alpha_o\cos\alpha \cosh\eta,\quad \quad C(\eta)=-1-\alpha_o^2 -\cos^2\alpha\sinh^2\eta  ,
\end{equation}
and  $\alpha_o=\xi_o/\xi$. The first term in the right-hand side of (\ref{eq:automur})
gives the constant Planck spectrum that is observed when the observer moves in a homogeneous medium. 
The second term gives the perturbation of the Planck spectrum that depends on $\tau$ and that is induced by the obstacle.
The calculation of the Wigner transform requires the following lemma.

\begin{lemma}
\label{lem:1}
For $a\in [0,1)$, $b \in \RR$, and $c\leq -1$, with $a+c+|b|<0$, we define the integral
\begin{equation}
\label{eq:lem1a}
\Psi(v;a,b,c) = \int_\RR \frac{\exp(iv s)}{a\cosh^2(s) +b \cosh(s)+c }ds .
\end{equation}
\begin{enumerate}
\item
If $a\in (0,1)$ and $b \neq 0$, then 
\begin{equation}
\label{eq:lem1b}
\Psi(v;a,b,c) = - \frac{2\pi}{\sqrt{\Delta}} \Big\{
\frac{ \sin \big[ v \, {\rm argcosh}(c_+) \big]}{ \tanh (\pi v) \sqrt{c_+^2-1}}
+
\frac{ \sin \big[ v \, {\rm argcosh}(|c_-|) \big]}{ \sinh (\pi v)\sqrt{c_-^2-1}}\Big\}  ,
\end{equation}
where
$$
c_\pm = \frac{-b\pm\sqrt{\Delta}}{2a}
,\quad\quad
\Delta = b^2-4ac .
$$
\item
If $a=0$ and $b \neq 0$, then 
\begin{equation}
\label{eq:lem1c}
\Psi(v;0,b,c) = - \frac{2\pi}{\sqrt{c^2-b^2} \sinh (\pi v)}  
  \sin \big[ v \, {\rm argcosh}(\sqrt{c/b}) \big] .
\end{equation}

\item
If  $a\in (0,1)$ and $b=0$, then
\begin{equation}
\label{eq:lem1d}
\Psi(v;a,0,c) = - \frac{\pi}{\sqrt{c^2+ca} \tanh (\pi v/2)}  
  \sin \big[ v \, {\rm argcosh}(\sqrt{-c/a}) \big] .
\end{equation}

\item
If $a=0$ and $b=0$, then 
\begin{equation}
\label{eq:lem1e}
\Psi(v;0,0,c) = \frac{2\pi}{c} \delta(v) .
\end{equation}
\end{enumerate}
\end{lemma}
This lemma can be proved by the residue theorem {(see \ref{proof:lemma1})}
and it can be applied to prove the following proposition.

\begin{proposition}
The Wigner transform $W(\tau,\omega)$
 of the recorded signal is a deformed version of the Planck spectrum $W_o(\omega)$ defined by (\ref{def:Wonu}):
 \begin{equation}
W(\tau,\omega)=
W_o(\omega) 
\big[  1 -  R(c_o\tau/\xi , \xi \omega/c_o) \big] .
\label{eq:Wmur}
\end{equation}
The correction $R$ has the following form:\\
- If the observer {trajectory} is not normally incident, $\alpha \neq 0$, then 
\begin{align}
\nonumber
R(\eta,\nu)=& \frac{\tanh(\pi \nu)}{\nu\tanh(2\pi\nu)}
\frac{ \sin \big[ 2\nu \, {\rm argcosh}(C_+(\eta)) \big] }{\sqrt{B^2(\eta)-4AC(\eta)} \sqrt{C_+^2(\eta)-1}}\\
&
+\frac{\tanh(\pi\nu)}{\nu\sinh(2\pi\nu)}
\frac{ \sin \big[ 2\nu \, {\rm argcosh}(|C_-(\eta)|) \big]}{\sqrt{B^2(\eta)-4AC(\eta)} \sqrt{C_-^2(\eta)-1}}  ,
\end{align} 
where
$$
C_\pm(\eta) = \frac{-B(\eta)\pm\sqrt{B^2(\eta)-4AC(\eta)}}{2A}   ,
$$
and $A,B(\eta),C(\eta)$ are given by (\ref{eq:defABC}).\\
- If the observer {trajectory} is normally incident, $\alpha = 0$,  
then
\begin{align}
R(\eta,\nu)=&
\left\{
\begin{array}{ll}
\displaystyle
\frac{\tanh(\pi\nu)}{\tanh(2\pi\nu)}
\frac{ \sin \big[ 2\nu \, {\rm argcosh}(|C(\eta)/B(\eta)|) \big]}{\nu \sqrt{C^2(\eta)-B^2(\eta)}}  
& \mbox{ if } \alpha_o \in (-1,0) , \\
\displaystyle
\frac{\tanh(\pi \nu)}{\sinh(2\pi\nu)}
\frac{ \sin \big[ \nu \, {\rm argcosh}(|C(\eta)/B(\eta)|) \big]}{2\nu \sqrt{C^2(\eta)-B^2(\eta)}}  
& \mbox{ if }  \alpha_o \in (0,+\infty) ,
 \\  0 & \mbox{ if } \alpha_o =0,
 \end{array}
 \right.
 \label{eq:expressR2}
\end{align} 
where $B(\eta),C(\eta)$ are given by (\ref{eq:defABC}).
\end{proposition}

Eq.~(\ref{eq:expressR2}) can also be written as
\begin{align}
R(\eta,\nu) =
\left\{
\begin{array}{ll}
\displaystyle
 \Big(1-\frac{1}{2\cosh^2(\pi\nu)}\Big) \frac{ \sin \big[ 2\nu \, \ln(\cosh(\eta)/|\alpha_o|) \big]}{ \nu \big(\cosh^2(\eta)-\alpha_o^2\big)}
& \mbox{ if } \alpha_o \in (-1,0) , \\
\displaystyle
\frac{1}{2\cosh^2(\pi\nu)}  \frac{ \sin \big[ 2\nu \, \ln(\alpha_o/\cosh(\eta)) \big]}{ \nu \big(\alpha_o^2-\cosh^2(\eta)\big)}
& \mbox{ if } \alpha_o \in (0,+\infty)
 , \\ 0 & \mbox{ if } \alpha_o =0
.
 \end{array}
 \right.
 \end{align}
This proposition shows how {the position of the obstacle} relative to the observer is encoded in the perturbation of the Wigner transform
of the recorded signal.

\subsection{The localization of the obstacle}
The obstacle can be detected by the observer except  in the special situation when 
$\alpha=0$ and $\xi_o=0$. We then have:
\begin{align}
\nonumber
\left< U(\tau+\frac{\tau'}{2})U(\tau-\frac{\tau'}{2})\right>
=
- \frac{c_o^2 f_o}{16\pi^2 \xi^2}\frac{1}{\sinh^2(\eta'/2)} 
- \frac{c_o^2 f_o}{16 \pi^2 \xi^2}
\frac{1}{\cosh^2(\eta)}  ,
\end{align}
where $\eta= c_o \tau /\xi$ and $\eta'=c_o\tau'/\xi$,
so that the perturbation of the Wigner transform
 is proportional to $\delta(\omega)$, which does not affect the spectrum at any $\omega \neq 0$.
 This was observed in \cite{rovelli}, but we show here that this result only hold true for a very particular situation,
 when the Rindler trajectory is normally incident ($\alpha=0$) and stops at exactly the distance $\xi$ from the 
 obstacle ($\xi_0=0$).

In general,
when $|\eta|\gg 1$, i.e. when the observer is far from the obstacle, then $C_\pm(\eta)\sim \cosh(\eta)$ and $R(\eta,\nu)$
is negligible.

\begin{figure}
\begin{center}
\begin{tabular}{c}
\includegraphics[width=5.4cm]{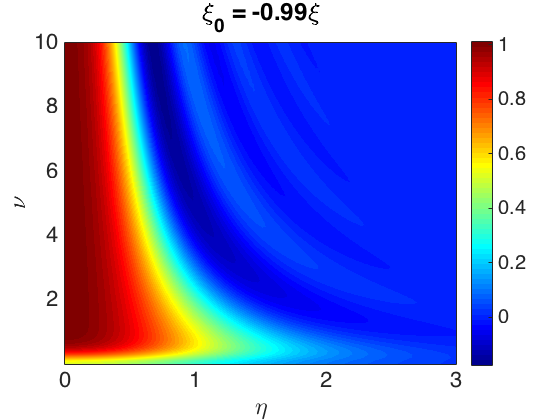}
\includegraphics[width=5.4cm]{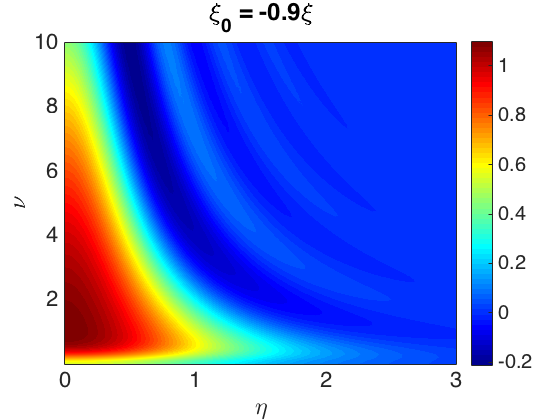}\\
\includegraphics[width=5.4cm]{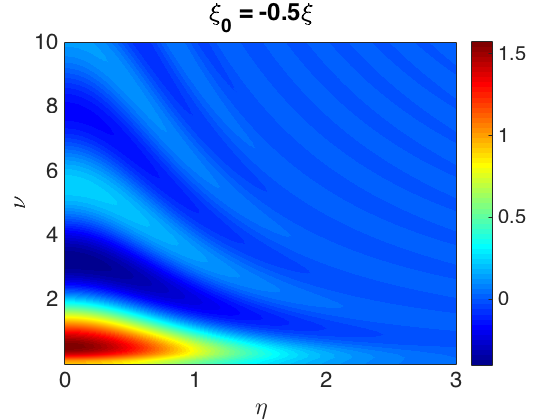}
\includegraphics[width=5.4cm]{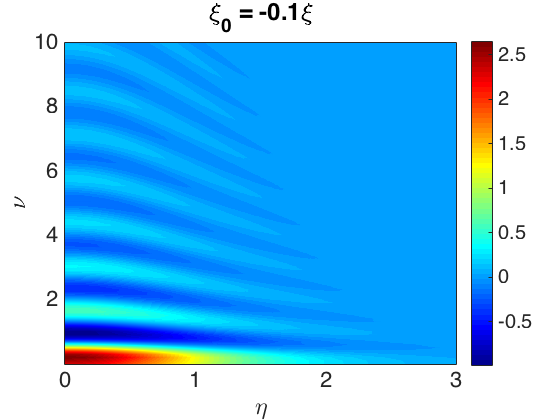}\\
\includegraphics[width=5.4cm]{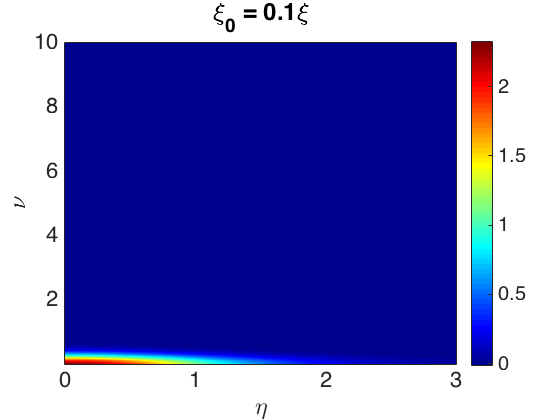}
\includegraphics[width=5.4cm]{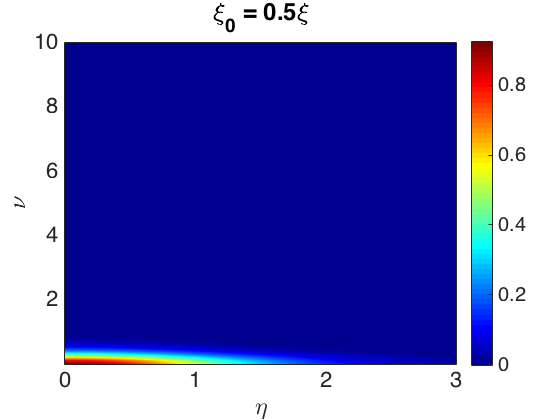}
\end{tabular}
\end{center}
\caption{
Correction $R(\eta,\nu)$ when $\alpha=0$ (normal incidence). 
{The correction $R$ is defined in (\ref{eq:Wmur}), as a function of $\eta=c_o \tau/\xi$ and $\nu=\xi \omega/c_o$.}
When $\xi_o=-0.99\xi$ and $\eta=0$, the observer is at distance $0.01\xi$ from the obstacle {(the trajectory is given by (\ref{eq:rindlerobs}))}.
\label{fig:rindler3}
}
\end{figure}

In Figure \ref{fig:rindler3} the correction $R(\eta,\nu)$ is plotted when the trajectory has normal incidence $\alpha=0$.
When the observer comes close to the obstacle (i.e. when $\xi_o $ is close to $-\xi$ and $\eta\simeq 0$)
the correction is close to one for a large band of frequencies:
$$
R(\eta,\nu) \stackrel{\xi_o \simeq -\xi, \, \eta\simeq 0}{\simeq} 1-\frac{1}{2\cosh^2(\pi\nu)} ,
$$ 
{which makes it easy for the observer to detect the correction, hence the obstacle.
The fact that the correction is approximately one comes from the Dirichlet boundary condition at the boundary of the obstacle, which makes
the field approximately zero close to the boundary.}
When the observer is far from the obstacle (for large $\eta$) the correction is close to zero.
When the observer is in the neighborhood of the obstacle, the correction has a complicated frequency-dependent structure.
The maximal correction is reached at positive frequency when $\xi_o<0$ and at zero-frequency when $\xi_o>0$.
The correction can be larger than one. There is no contradiction, the Wigner transform can take (locally in $\eta$) negative values.

\begin{figure}
\begin{center}
\begin{tabular}{c}
\includegraphics[width=5.4cm]{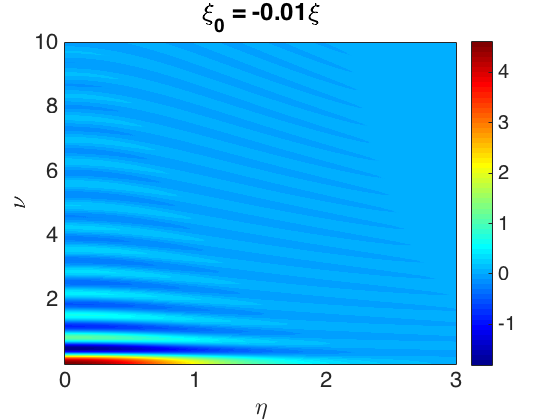}
\includegraphics[width=5.4cm]{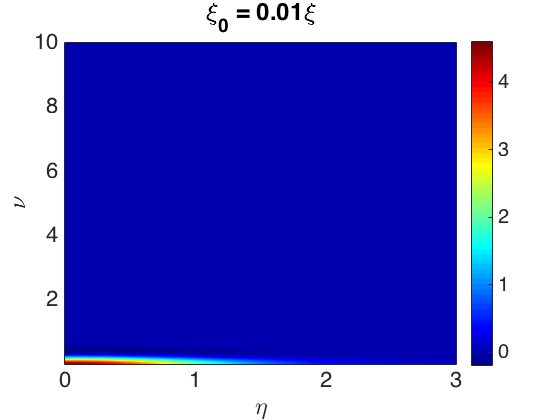}\\
\includegraphics[width=5.4cm]{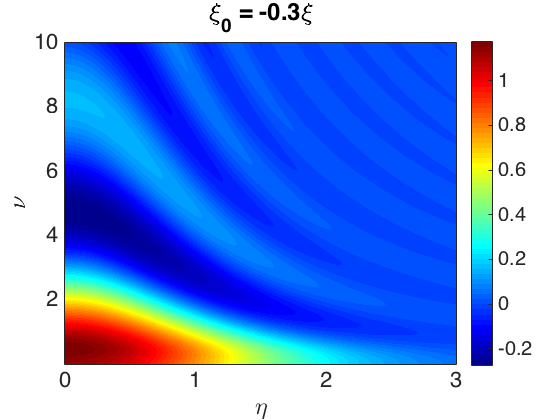}
\includegraphics[width=5.4cm]{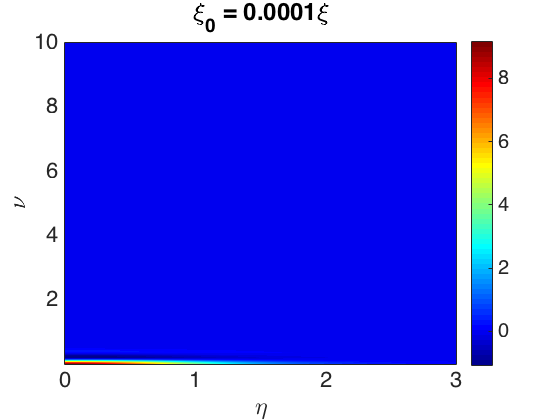}
\end{tabular}
\end{center}
\caption{Correction $R(\eta,\nu)$ when $\alpha=0$  (normal incidence) and $\alpha_o = \xi_o/\xi$ is close to zero.
The correction becomes proportional to $\delta(\nu)$.
    \label{fig:rindler3b}
}
\end{figure}

In Figure \ref{fig:rindler3b}, the correction $R(\eta,\nu)$ is plotted when the trajectory has normal incidence $\alpha=0$
and $\alpha_o$ is very close to zero. This illustrates the above remark that the correction is in this case concentrated
at very small frequencies, and that in the limit $\alpha_o\to 0$ 
{it is} proportional to $\delta(\nu)$.

\begin{figure}
\begin{center}
\begin{tabular}{c}
\includegraphics[width=5.4cm]{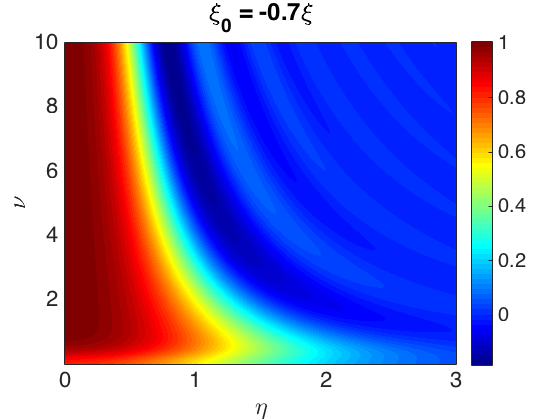}
\includegraphics[width=5.4cm]{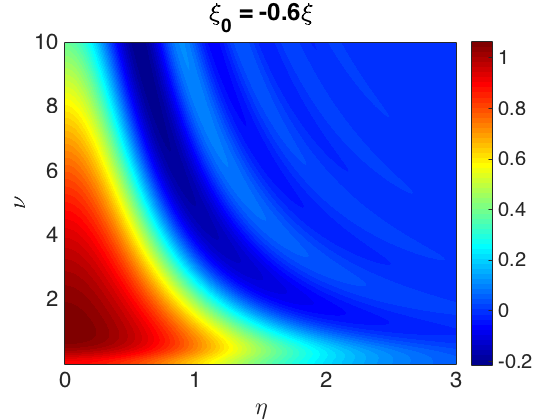}\\
\includegraphics[width=5.4cm]{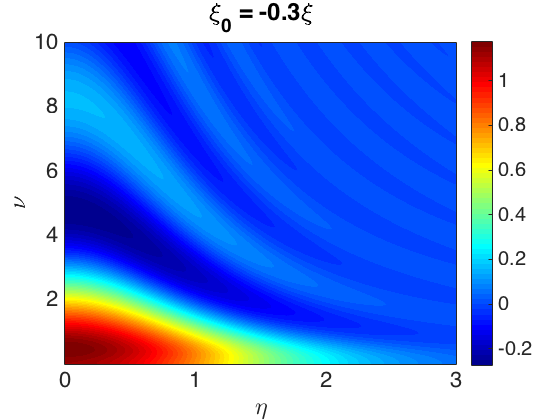}
\includegraphics[width=5.4cm]{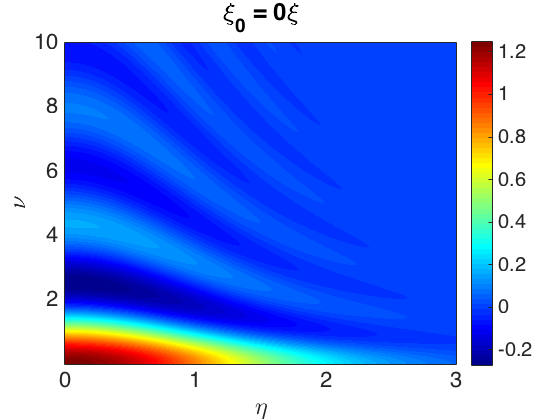}\\
\includegraphics[width=5.4cm]{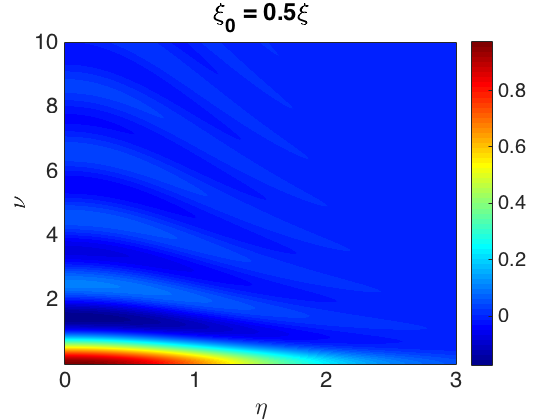}
\includegraphics[width=5.4cm]{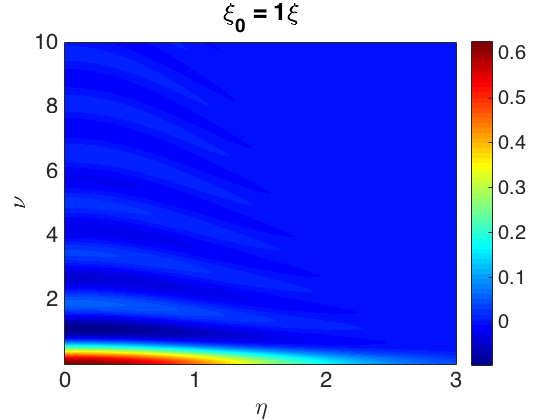}
\end{tabular}
\end{center}
\caption{Correction $R(\eta,\nu)$ when $\alpha=\pi/4$  (oblique incidence).
When $\xi_o=-0.7\xi$ and $\eta=0$, the observer is at distance $0.0071\xi$ from the obstacle.
    \label{fig:rindler4}
}
\end{figure}

In Figure \ref{fig:rindler4} the correction $R(\eta,\nu)$ is plotted when the trajectory has oblique incidence $\alpha=\pi/4$.
The results are quantitatively different, but qualitatively similar. In particular, when the observer is very close to the 
obstacle ($\xi_o \simeq - \xi \cos(\alpha)$, $\eta \simeq 0$) then the correction is close to one for a large band of frequencies
because of the Dirichlet condition:
$$
R(\eta,\nu) \stackrel{\xi_o \simeq -\xi \cos \alpha, \, \eta\simeq 0}{\simeq} 1-\frac{1}{2\cosh^2(\pi\nu)}
\Big( 1 - \frac{\sin \big[2\nu {\rm argcosh}(1+2/\tan^2\alpha)\big]}{4\nu \sqrt{1/\tan^2\alpha +1 /\tan^4\alpha}}\Big)  .
$$

In the point of view of the inverse problem, one could implement an optimal control strategy 
to minimize the least square mismatch between
an observed spectrum correction and the theoretical spectrum correction parameterized by $\alpha$ and $\alpha_o$.
This can be done for a given time $\tau$, or for several $\tau$ in order to improve the 
resolution and the robustness of the estimation method.
These are standard procedures \cite{ammari}.

\section{Conclusion}
The results reported in this paper
show that, as long as a Rindler observer is far from an obstacle,
it perceives a constant Planck spectrum when the illumination comes from noise sources with Lorentz-invariant spectrum.
However, when the observer comes into the neighborhood of an obstacle,
the Wigner transform of the recorded signal deviates from the Planck spectrum
and the deviation can be used to localize the obstacle.
This demonstrates that a passive observer can perceive its environment from the signal that it records
and that is transmitted by  noise sources.

{In this paper the obstacle has the form of an infinite perfect mirror.
It would be interesting to address more general obstacles, for which 
the curvature of the obstacle surface would be important
\cite{deutch}.}
This result could also be extended qualitatively to other trajectories: 
an obstacle would modify the Wigner transform of the recorded signal.
However the Rindler trajectory is the unique straight-line trajectory  that  allows the observer
to feel a constant spectrum whatever its time or position, as long as the observer is far from any obstacle.
It is only when it comes into the neighborhood of an obstacle that the spectrum is modified,
which allows the observer to detect the obstacle and to determine its relative position with respect to the obstacle
 {once the entire trajectory has been traversed.}
Finally, this paper only addresses noise sources with Lorentz-invariant spectra.
If we look for other applications, other spectra may be more appropriate.
 {Remarks \ref{rem:1} and \ref{rem:2} seem to indicate that our main results 
are somewhat robust with respect to the form of the source spectrum and the duration of the recording time window,
but more detailed work is needed to quantify the loss of accuracy and resolution in these  general cases. }

\section*{Acknowledgments}
This work was partially supported by LABEX WIFI (Laboratory of Excellence ANR-10-LABX-24) 
within the French Program Investments for the Future under reference ANR-10-IDEX-0001-02 PSL*.

\appendix
\section{Helmholtz-Kirchhoff identity}
For any $\bx_1,\bx_2 \in \RR^3$ we have for  $L \to \infty$:
\begin{equation}
\label{hk}
{\rm Im}
 \hat{G}(\omega, \bx_1,\bx_2)  = \frac{ \omega}{c_o}   \int_{\partial B({\bf 0},L)} 
\overline{ \hat{G}( \omega,\bx_1,\by)}  \hat{G}( \omega ,\bx_2,\by)  d \sigma(\by) 
 ,
\end{equation}
where $d \sigma(\by) $ is the surface integral.
It is a consequence of second Green's identity and Sommerfeld radiation condition \cite{garpapa16}.

\section{{Proofs}}
\subsection{Proof of (\ref{eq:expresauto3dot})}
\label{app:B1}
We consider (\ref{eq:expresauto3do}) with $\hat{F}(\omega)=|\omega|$ and $\tau'\neq 0$.
We denote $\eta= c_o \tau /\xi$ and $\eta'=c_o\tau'/\xi$.
By (\ref{eq:rindler3d}), we have 
\begin{align*}
 \big|\bX(\tau+\frac{\tau'}{2})- \bX(\tau-\frac{\tau'}{2})\big| =& 2 \xi \big| \sinh( \eta)\sinh(\frac{\eta'}{2})\big| , \\
 T(\tau+\frac{\tau'}{2})-T(\tau-\frac{\tau'}{2})\big) =&  2\xi  \cosh( \eta)\sinh(\frac{\eta'}{2}) ,
\end{align*}
so that
 \begin{align*}
\left< U(\tau+\frac{\tau'}{2})U(\tau-\frac{\tau'}{2})\right>
=&
\frac{1}{4\pi^2}  \frac{1}{\frac{2\xi}{c_o}  \sinh( \eta)\sinh(\frac{\eta'}{2}) } \\
& \times \int_0^\infty
\sin \Big( \frac{\omega}{c_o} 2 \xi  \sinh( \eta)\sinh(\frac{\eta'}{2}) \Big)
\cos \Big(  \frac{\omega}{c_o} 2 \xi  \cosh( \eta)\sinh(\frac{\eta'}{2}) \Big) d\omega \\
=&
\frac{1}{8\pi^2}  \frac{1}{\frac{2\xi}{c_o}  \sinh( \eta)\sinh(\frac{\eta'}{2}) }  \int_0^\infty
\sin \Big( \frac{\omega}{c_o} 2 \xi  \sinh(\frac{\eta'}{2}) [\sinh( \eta)+\cosh( \eta)]\Big)\\
&
+
\sin \Big( \frac{\omega}{c_o} 2 \xi  \sinh(\frac{\eta'}{2}) [\sinh( \eta)-\cosh( \eta)]\Big) d\omega.
\end{align*} 
We have {(for the justification of the inversion of the integral and the limit $\epsilon \to 0$, see \ref{app:B3})}
\begin{align}
\label{eq:inversion}
 \int_0^\infty \sin (\omega t) d\omega =  {\rm Im} \lim_{\epsilon \to 0^+} \int_0^\infty e^{(i t -\epsilon) \omega} d\omega = {\rm Im} \frac{-1}{i t -\epsilon} = \frac{1}{t},
 \end{align}
 so that
  \begin{align*}
\left< U(\tau+\frac{\tau'}{2})U(\tau-\frac{\tau'}{2})\right>
=&
\frac{1}{8\pi^2}  \frac{1}{\frac{2\xi}{c_o}  \sinh( \eta)\sinh(\frac{\eta'}{2}) } \\
& \times \Big[ \frac{1}{ \frac{2\xi}{c_o} \sinh(\frac{\eta'}{2}) [\sinh( \eta)+\cosh( \eta)]}
+
\frac{1}{ \frac{2\xi}{c_o} \sinh(\frac{\eta'}{2}) [\sinh( \eta)-\cosh( \eta)]}\Big]  \\
=&
- \frac{c_o^2}{16 \pi^2 \xi^2 \sinh^2 (\frac{\eta'}{2})}.
\end{align*}

\subsection{Proof of (\ref{def:Wonu0})}
\label{app:B2}
We show that the inverse Fourier transform of $\omega \mapsto W(\tau,\omega)$ defined by (\ref{def:Wonu0}) 
gives (\ref{eq:expresauto3dot}) for any $\tau'\neq 0$ (the proof follows closely \cite{boyer80}):
\begin{align*}
\frac{1}{2\pi }\int_\RR W(\tau,\omega) e^{- i\omega \tau'}d\omega
&= 
\frac{f_o}{4\pi^2} \int_0^\infty  \frac{\omega}{\tanh(\pi \xi \omega /c_o)} \cos(\omega \tau') d\omega \\
&= 
\frac{f_o}{4\pi^2}
\Big\{
 \int_0^\infty  \omega \cos(\omega \tau') d\omega 
+ \int_0^\infty  \frac{2 \omega \cos (\omega \tau')}{\exp(2\pi \xi \omega /c_o) -1}d\omega \Big\} .
\end{align*}
The singular part can be computed by
 \begin{equation}
 \label{eq:inversion2}
 \int_0^\infty \omega \cos (\omega t) d\omega =  {\rm Re} \lim_{\epsilon \to 0^+} \int_0^\infty \omega e^{(i t -\epsilon) \omega} d\omega = -\frac{1}{t^2},
\end{equation}
and by using \cite[formula 3.951.5]{grad} we get the expression of the regular part, so that we obtain
\begin{align*}
\frac{1}{2\pi }\int_\RR W(\tau,\omega) e^{- i\omega \tau'}d\omega
&= 
\frac{f_o}{4\pi^2}
\Big\{
- \frac{1}{{\tau'}^2}
+  \Big[ \frac{1}{{\tau'}^2}  - \frac{\frac{c_o^2}{4\xi^2}}{\sinh^2 ( c_o \tau'/(2\xi))} \Big]  \Big\},
\end{align*}
 which gives (\ref{def:Wonu0}).

{
\subsection{Finite-energy spectrum}
\label{app:B3}
Here we revisit the two previous appendices when the source spectrum
is of the form 
\begin{equation}
\label{eq:specepsilon}
\hat{F}(\omega) = f_o |\omega|\exp(-\epsilon |\omega|)
\end{equation}
for some $\epsilon>0$.
The goal is twofold:
First we want to justify the inversion of the integral and the limit $\epsilon \to 0$ in (\ref{eq:inversion}) and (\ref{eq:inversion2}).
Second we want to show that we can deal with a noise source spectrum with finite energy 
and amplitude, and therefore classical recorded signals with finite energy and amplitude, without altering the results.
}

{
By (\ref{eq:rindler3d}), we have 
 \begin{align*}
\left< U(\tau+\frac{\tau'}{2})U(\tau-\frac{\tau'}{2})\right>
=&
\frac{1}{8\pi^2}  \frac{1}{\frac{2\xi}{c_o}  \sinh( \eta)\sinh(\frac{\eta'}{2}) }  \int_0^\infty
\Big[ \sin \Big( \frac{\omega}{c_o} 2 \xi  \sinh(\frac{\eta'}{2}) e^{\eta}\Big)\\
&
-
\sin \Big( \frac{\omega}{c_o} 2 \xi  \sinh(\frac{\eta'}{2}) e^{-\eta} \Big) \Big]\exp(- \epsilon \omega)  d\omega.
\end{align*} 
Then, using 
$$
 \int_0^\infty \sin (\omega t) \exp(-\epsilon \omega) d\omega =  {\rm Im} \int_0^\infty e^{(i t -\epsilon) \omega} d\omega = {\rm Im} \frac{-1}{i t -\epsilon} = \frac{t}{t^2+\epsilon^2},
$$
we find that, for any $\epsilon>0$:
  \begin{align}
\left< U(\tau+\frac{\tau'}{2})U(\tau-\frac{\tau'}{2})\right>
=&
\frac{1}{8\pi^2 \sinh( \eta)}    \Big[ \frac{ e^\eta}{ \big(\frac{2\xi}{c_o} \big)^2\sinh(\frac{\eta'}{2})^2 e^{2\eta} +\epsilon^2}
 -
\frac{ e^{-\eta}}{ \big(\frac{2\xi}{c_o} \big)^2 \sinh(\frac{\eta'}{2})^2 e^{-2\eta} +\epsilon^2}\Big]  .
\label{eq:Ueps}
\end{align}
The limit $\epsilon \to 0$ can be taken in this expression to get 
that (\ref{eq:expresauto3dot}) gives the correct value of the autocorrelation function of the recorded signal
when the source spectrum is (\ref{eq:specepsilon}) and $\sinh(c_o  |\tau'| /(2\xi) ) >  (c_o\epsilon/\xi) \exp( c_o |\tau| /\xi)  $,
i.e. $|\tau'|>O(\epsilon)$.
}

{
A similar analysis can be carried out for the Wigner transform of the recorded signal (\ref{def:wigner}).
Using \cite[formula 3.983.1]{grad}, we get that, for any $\epsilon >0$,
\begin{align}
\nonumber
W(\tau,\omega) = &
\frac{c_o}{4\pi \xi (e^{2\eta}-1)} 
\frac{
\sinh\big( \frac{\omega \xi}{c_o}{\rm arccos}\big(-1+\frac{c_o^2}{2 \xi^2} e^{-2\eta} \epsilon^2\big)\big)
}
{
\sqrt{ 1-(1-\frac{c_o^2}{2 \xi^2} e^{-2\eta} \epsilon^2)^2} \sinh \big( \frac{\omega \xi \pi}{c_o}\big)
}
\\
&+
\frac{c_o}{4\pi \xi (e^{-2\eta}-1)} \frac{
\sinh\big( \frac{\omega \xi}{c_o}{\rm arccos}\big(-1+\frac{c_o^2}{2 \xi^2} e^{2\eta} \epsilon^2\big)\big)
}{
\sqrt{ 1-(1-\frac{c_o^2}{2 \xi^2} e^{2\eta} \epsilon^2)^2} \sinh \big( \frac{\omega \xi \pi}{c_o}\big)
}
.
\label{eq:Weps}
\end{align}
The limit $\epsilon \to 0$ can be taken in this expression (using ${\rm arccos}(-1+s) = \pi- \sqrt{2s}+O(s^{3/2})$ and $\sqrt{1-(1-s)^2} = \sqrt{2s}+O(s^{3/2})$ as $s\to0$) 
to get 
that (\ref{def:Wonu0}) gives the correct value of the Wigner transform of the recorded signal
when the source spectrum is (\ref{eq:specepsilon}) and
$c_o/\xi+|\omega|    <\epsilon^{-1} \exp(- c_o |\tau| /\xi) $,
i.e. $|\omega|<O(\epsilon^{-1})$.
}

\subsection{Proof of Proposition \ref{prop:2}}
\label{proof:prop2}%
A stationary observer records the signal $U(\tau)=u(\tau,\bx_0)$ whose autocorrelation function is
$$
\left< U(\tau+\frac{\tau'}{2})U(\tau-\frac{\tau'}{2})\right>=
\frac{1}{8\pi^2}
\int_\RR \hat{F}(\omega) \exp(i \omega \tau) d\omega .
$$
Let us consider an observer moving along the $z$-axis at the constant velocity $v >0$.
The trajectory with proper time $\tau$ is of the form
$(t(\tau),{\bx}(\tau))$ with ${\bx}(\tau)=(0,0,z(\tau))$. The time in the laboratory frame $t$ is related to 
the proper time $\tau$ by
\begin{equation}
\label{eq:mink}
\dot{t}^2 - \dot{z}^2 /c_o^2=1,
\end{equation}
with
$\dot{t}=\partial_\tau t$, $\dot{z}=\partial_\tau z$ (this comes from the fact that the Minkowski metric 
is $ds^2 = c_o^2dt^2 - dz^2$ and the proper time $\tau$ is $s/c_o$).
Without loss of generality (since the medium is invariant by any spatial shift), we can assume $t(0)=0$ and $z(0)=0$.
As the velocity $v$ is constant, we have $\frac{dz}{dt}=v$, or $\dot{z}=v\dot{t}$. 
Substituting into (\ref{eq:mink}) this gives $\dot{t}=\gamma$, with $\gamma=1/\sqrt{1-v^2/c_o^2}$ (Lorentz factor), and therefore
$t(\tau)=\gamma \tau$ and $z(\tau)=\gamma v \tau$.
From (\ref{eq:exprescov3do}) the autocorrelation function of the recorded signal $U(\tau)=u(t(\tau),{\bx}(\tau))$ has the form
\begin{align}
\left< U(\tau+\frac{\tau'}{2})U(\tau-\frac{\tau'}{2})\right>
=
\frac{1}{8\pi^2} \int_\RR  \hat{F}(\omega) 
{\rm sinc}\Big(\frac{\omega}{c_o} \gamma v \tau' \Big) \exp(i \omega \gamma \tau') d\omega ,
\end{align}
which is a function of $\tau'$ only, which may depend on $v$.
We want to identify the source spectrum $\hat{F}$ such that the autocorrelation function (or equivalently the Wigner transform)
does not depend on $v$.
Using the fact that, for any $\alpha> 0$, 
$$
\int_\RR {\rm sinc}(\alpha s) \exp(-i\Omega s) ds= \frac{\pi}{\alpha} {\bf 1}_{[-\alpha,\alpha]} (\Omega),
$$
we find that the Wigner transform is, for $\omega>0$,
\begin{equation}
\label{eq:expressWv}
W(\tau,\omega) = \frac{1}{8\pi \gamma \frac{v}{c_o}} \int_{\frac{\omega}{\gamma(1+\frac{v}{c_o})}}^{\frac{\omega}{\gamma(1-\frac{v}{c_o})}}
\hat{F}_1(\omega') d\omega' ,
\end{equation}
with $\hat{F}_1(\omega) = \hat{F}(\omega)/\omega$.
The Wigner transform can be expanded for small $v/c_o$ as
$$
W(\tau,\omega) = \frac{1}{4\pi} \hat{F}_1(\omega) + \frac{v^2}{8\pi c_o^2}\Big[\omega \partial_\omega \hat{F}_1(\omega) +\frac{\omega^2}{3} \partial_\omega^2  \hat{F}_1(\omega) \Big]+o\Big(\frac{v^2}{c_o^2}\Big).
$$
Therefore, 
a necessary condition for the Wigner transform to be independent of $v$  is that $\hat{F}_1$ should satisfy $\omega \partial_\omega \hat{F}_1(\omega) +\frac{\omega^2}{3} \partial_\omega^2  \hat{F}_1(\omega)=0$, that is to say, $\hat{F}_1(\omega)$ should be of the form
$ \hat{F}_1(\omega) = f_o+ \frac{f_1}{\omega^2}$, or
\begin{equation}
\label{eq:formhatFLorentz}
\hat{F}(\omega) = f_o |\omega|+ \frac{f_1}{|\omega|} .
\end{equation}
It turns out that (\ref{eq:formhatFLorentz}) is also a sufficient condition.
Indeed, if $\hat{F}(\omega)=f_o |\omega|$, then (\ref{eq:expressWv}) gives (for $\omega>0$)
$$
W(\tau,\omega) =\frac{f_o }{8\pi  \gamma \frac{v}{c_o}}
\Big( \frac{\omega}{\gamma(1-\frac{v}{c_o})} - \frac{\omega}{\gamma(1+\frac{v}{c_o})}
\Big) = \frac{f_o \omega}{4\pi} ,
$$
and if $\hat{F}(\omega)=f_1/|\omega|$, then (\ref{eq:expressWv}) gives (for $\omega>0$)
$$
W(\tau,\omega) =\frac{f_1}{8\pi \gamma \frac{v}{c_o}}
\Big( \frac{\gamma(1+\frac{v}{c_o})} {\omega}- \frac{\gamma(1-\frac{v}{c_o})} {\omega} \Big) = \frac{f_1}{4\pi \omega} .
$$
This completes the proof of Proposition \ref{prop:2}.

 \subsection{Proof of Proposition \ref{prop:3a}}
 \label{proof:prop3a}
We have
$$
\left< U(\tau+\frac{\tau'}{2})U(\tau-\frac{\tau'}{2})\right>=
\int_{\RR^3} {\cal A}(|\bk|)  \exp \big[ i 2\xi \sinh(\eta'/2) \big( |\bk| \cosh (\eta) - k_z \sinh(\eta)\big) \big] d\bk ,
$$
with $\eta= c_o \tau /\xi$ and $\eta'=c_o\tau'/\xi$.
After the change of variable $\bk' = (k_x,k_y,k_z \cosh(\eta) -|\bk| \sinh(\eta))$, we
get
$$
\left< U(\tau+\frac{\tau'}{2})U(\tau-\frac{\tau'}{2})\right>=
\int_{\RR^3}  \frac{{\cal A}( {\cal K}(\bk',\eta)){\cal K}(\bk',\eta)}{|\bk'|} \exp \big[ i 2\xi \sinh(\eta'/2)  |\bk'|   \big] d\bk' ,
$$
with 
$$
{\cal K}(\bk',\eta) =  \cosh(\eta) |\bk'|+ k_z' \sinh(\eta) .
$$
This function does not depend on $\eta$ if and only if $ k \mapsto {\cal A}(k) k$ is constant,  if and only if  $\hat{F}(\omega)$
is proportional to $|\omega|$ by (\ref{eq:calAk}).
This completes the proof of Proposition \ref{prop:3a}.

\subsection{Proof of Lemma \ref{lem:1}}
\label{proof:lemma1}
We address case 1 in order to compute $\Psi(v;a,b,c)$ defined by  (\ref{eq:lem1a}).
We have 
$$
\frac{1}{a\cosh^2(s) +b \cosh(s)+c }
=\frac{1}{a}  \frac{1}{\cosh (s) -c_+} \frac{1}{\cosh(s) -c_-}  ,
$$
with
$$
c_\pm = \frac{-b \pm \sqrt{\Delta}}{2a} , \quad \Delta=b^2-4ac.
$$
We can check that $c_+>1$ (because $-c-a-b>0$) and $c_- <-1$  (because $-c-a+ b>0$).
We denote
$$
x_\pm = \pm {\rm argcosh}(c_+),
\quad
\tilde{x}_\pm = i\pi  \pm {\rm argcosh}(|c_-|),
\quad
\check{x}_\pm  = 2 i \pi  + x_\pm.
$$
We have $\cosh(x_\pm)= \cosh(\check{x}_\pm)=c_+$ and $\cosh(\tilde{x}_\pm)=-c_-$.
We apply the residue theorem on the closed contour which is a rectangle 
$[-M,M] \cup [M,M+2i\pi]\cup[M+2i\pi,-M+2i\pi]\cup[-M+2i\pi,-M]$, with $M\to +\infty$. The contour contains two poles   $\tilde{x}_\pm$
 and it passes through four poles $x_\pm$ and $\check{x}_\pm$,
so we get by the residue theorem:
\begin{align*}
\Psi(v;a,b,c) (1-e^{-2\pi v}) =&  \frac{2i\pi}{a (c_+-c_-)} 
\Big\{
\frac{1}{2} \frac{e^{i v x_-}}{\sinh(x_-)}
+
\frac{1}{2} \frac{e^{i v x_+}}{\sinh(x_+)}
-
\frac{e^{i v \tilde{x}_+}}{\sinh(\tilde{x}_+)} \\
&
-
\frac{e^{i v \tilde{x}_-}}{\sinh(\tilde{x}_-)}
+
\frac{1}{2}
\frac{e^{i v \check{x}_+ }}{ \sinh( \check{x}_+)}
+
\frac{1}{2}
\frac{e^{i v  \check{x}_- }}{  \sinh(\check{x}_-)}
\Big\} ,
\end{align*}
which gives the desired result  (\ref{eq:lem1b}) using $\sinh(x_\pm) = \pm \sqrt{c_+^2-1}$,
$\sinh(\tilde{x}_\pm) = \mp \sqrt{c_-^2-1}$, $\sinh(\check{x}_\pm)= \pm  \sqrt{c_+^2-1}$, and $a (c_+-c_-) = \sqrt{\Delta}$.

Next we address case 2. We have
$$
\Psi(v;0,b,c) =  \frac{1}{b} \int_\RR \frac{\exp(i v s)}{\cosh(s)+c'} ds, 
$$
with $c'=c/b >1$.
We denote $x_\pm = i \pi +{\rm argcosh}(c')$. We have $\cosh(x_\pm )+c'=0$.
We apply the residue theorem on the contour $[-M,M] \cup [M,M+2i\pi]\cup[M+2i\pi,-M+2i\pi]\cup[-M+2i\pi,-M]$, with $M\to +\infty$, which contains the two poles $x_\pm$,
and we get
\begin{align*}
&
\Psi(v;0,b,c)  (1-e^{-2\pi v}) =  \frac{2i\pi}{b} 
\Big\{
\frac{e^{i v x_+}}{  \sinh(x_+)}
+
\frac{e^{i v x_-}}{  \sinh(x_-)} \Big\} ,
\end{align*}
which gives the desired result  (\ref{eq:lem1c}).

Finally we address case 3. We have
$$
\Psi(v;a,0,c) = \frac{1}{a} \int_\RR \frac{\exp( i (v/2) s)}{\cosh( s) +c'} ds
$$
with $c'= 1+2c/a< -1$.
We denote $x_\pm = \pm {\rm argcosh}(|c'|)$ and $\tilde{x}_\pm = 2i\pi +x_\pm$. 
We have $\cosh(x_\pm )+c'=\cosh(\tilde{x}_\pm )+c'=0$.
We apply the residue theorem on the contour $[-M,M] \cup [M,M+2i\pi]\cup[M+2i\pi,-M+2i\pi]\cup[-M+2i\pi,-M]$, with $M\to +\infty$, which passes through the four poles $x_\pm$ and $\tilde{x}_\pm$, 
and we get
\begin{align*}
&
\Psi(v;a,0,c) (1-e^{-\pi v}) =  \frac{i\pi}{a} 
\Big\{
\frac{e^{i (v/2)x_+}}{  \sinh(x_+)}
+
\frac{e^{i (v/2) x_-}}{ \sinh(x_-)} 
+
\frac{e^{i (v/2)  \tilde{x}_+}}{  \sinh( {x}_+)}
+
\frac{e^{i (v/2) \tilde{x}_-}}{ \sinh({x}_-)}
\Big\} ,
\end{align*}
which gives the desired result  (\ref{eq:lem1d}) with $ {\rm argcosh}(|c'|)) = 2 {\rm argcosh}(\sqrt{-c/a})$.



\begin{thebibliography}{99}

\bibitem{ammari}
H. Ammari, J. Garnier, W. Jing, H. Kang, M. Lim, K. S\o lna, and H. Wang, 
Mathematical and Statistical Methods for Multistatic Imaging, 
Lecture Notes in Mathematics, Vol. 2098, Springer, Berlin, 2013.

\bibitem{antoniou12} 
M. Antoniou, Z. Zeng, L. Feifeng, and M. Cherniakov, 
Experimental demonstration of passive BSAR imaging
using navigation satellites and a fixed receiver, 
IEEE Geoscience and Remote Sensing Letters {\bf 9}, 477--481 (2012). 

\bibitem{boyer80}
T. H. Boyer,
Thermal effects of acceleration through random classical radiation,
Phys. Rev. D {\bf 21},  2137--2148 (1980).

\bibitem{boyer84}
T. H. Boyer,
Thermal effects of acceleration for a classical dipole oscillator in classical electromagnetic zero-point radiation,
Phys. Rev. D {\bf 29},  1089--1094 (1984).

\bibitem{cheney01}
M. Cheney,
A mathematical tutorial on synthetic aperture radar,
{SIAM Review} {\bf 43}, 301--312 (2001). 

\bibitem{crispino}
 L. C. B. Crispino, A. Higuchi, and G. E. A. Matsas, 
 The Unruh effect and its applications, 
 Reviews of Modern Physics {\bf 80}, 787--838 (2008).

\bibitem{curlander}
J. C. Curlander and R. N. McDonough, 
{\it Synthetic aperture radar}, Wiley, New York, 1991.

\bibitem{curtis06}
A. Curtis, P. Gerstoft, H. Sato, R. Snieder, and K. Wapenaar, 
Seismic interferometry turning noise into signal,
The Leading Edge {\bf 25}, 1082--1092 (2006). 

\bibitem{deutch}
{
D. Deutsch and P. Candelas,
Boundary effects in quantum field theory, 
Phys. Rev. D {\bf 20}, 3063--3080 (1979).
}

\bibitem{fink17}
M. Fink and J. Garnier, 
{Ambient noise correlation-based imaging with moving sensors}, 
Inverse Problems and Imaging {\bf 11}, 477--500 (2017). 

\bibitem{finkgar15}
J. Garnier and M. Fink, 
Super-resolution in time-reversal focusing on a moving source, 
Wave Motion {\bf 53}, 80--93 (2015).

\bibitem{garpapa09}
J. Garnier and G. Papanicolaou, 
Passive sensor imaging using cross correlations of noisy signals in a scattering medium, 
SIAM J. Imaging Sciences {\bf 2},  396--437 (2009).

\bibitem{garpapa15}
J. Garnier and G. Papanicolaou, 
Passive synthetic aperture imaging, 
SIAM J. Imaging Sciences {\bf 8}, 2683--2705  (2015). 

\bibitem{garpapa16}
J. Garnier and G. Papanicolaou,
\textit{Passive Imaging with Ambient Noise},
Cambridge University Press, Cambridge, 2016.

\bibitem{gouedard08}
P. Gou\'edard, L. Stehly, F. Brenguier, M. Campillo, Y. Colin de Verdi\`ere, E. Larose,  L. Margerin, P. Roux, F. J. Sanchez-Sesma, N. M. Shapiro, and R. L. Weaver,
{Cross-correlation of random fields: mathematical approach and applications},
{Geophysical Prospecting} {\bf 56}, 375--393 (2008). 

\bibitem{grad}
{ I. S. Gradshteyn et I. M. Ryzhik,
{\it Table of Integrals, Series, and Products},
Academic Press, San Diego, 1980.
}

\bibitem{book:ulf}
 U. Leonhardt, 
 \textit{Essential Quantum Optics: From Quantum Measurements to Black Holes}, 
 Cambridge University Press, Cambridge, 2010.

\bibitem{ulf}
U. Leonhardt, I. Griniasty, S. Wildeman, E. Fort, and M. Fink,
Classical analog of the Unruh effect,
Phys. Rev. A {\bf 98}, 022118 (2018). 

\bibitem{rindler}
W. Rindler, 
Kruskal space and the uniformly accelerated frame,
Am. J. Phys. {\bf 34}, 1174--1178 (1966).

 \bibitem{rodriguez10}
M. Rodriguez-Cassola, S. V. Baumgartner, G. Krieger, and A. Moreira, 
Bistatic TerraSAR-X/F-SAR Spaceborne-Airborne SAR experiment: Description, data processing, and results,
IEEE Transactions on Geoscience and Remote Sensing {\bf 48}, 781--794 (2010).

\bibitem{rovelli}
C. Rovelli and M. Smerlak, 
Unruh effect without trans-horizon entanglement, 
Phys. Rev. D {\bf 85}, 124055 (2012).

\bibitem{rytov}
{
S. M. Rytov,
Theory of the Electric Fluctuations and Thermal Radiation [in Russian], Publication of Acad. of Sciences of USSR, Moscow (1953), English translation: Air Force Cambridge Research Center, Bedford, MA (1959).
}

\bibitem{shapiro05}
N. M. Shapiro, M. Campillo, L. Stehly,  and M. H. Ritzwoller,
{High-resolution surface wave tomography from ambient noise},
{Science} {\bf 307}, 1615--1618 (2005). 

\bibitem{unruh}
W. G. Unruh, Notes on black-hole evaporation,
Phys. Rev. D {\bf 14}, 870--892 (1976).

\bibitem{wap10a}
{K. Wapenaar, D. Draganov, R. Snieder, X. Campman, and A. Verdel},
{Tutorial on seismic interferometry:
Part 1 - Basic principles and applications},
{Geophysics} {\bf 75}, A195--A209 (2010). 


\end{thebibliography}
\end{document}